\documentclass[reqno,11pt]{amsart}
\usepackage{amssymb, amsmath,latexsym,amsfonts,amsbsy, amsthm}
\usepackage{color}
\usepackage{xcolor}
\usepackage{mathrsfs}
\usepackage{cleveref}

\allowdisplaybreaks[4]

\let\Om=\Omega

\newcommand{\beq}{\begin{equation}}
	\newcommand{\eeq}{\end{equation}}
\newcommand{\ben}{\begin{eqnarray}}
	\newcommand{\een}{\end{eqnarray}}
\newcommand{\beno}{\begin{eqnarray*}}
	\newcommand{\eeno}{\end{eqnarray*}}

\newtheorem{remark}{Remark}[section]

\numberwithin{equation}{section}
\newtheorem{theorem}{Theorem}[section]
\newtheorem{definition}[theorem]{Definition}
\newtheorem{lemma}[theorem]{Lemma}
\newtheorem{proposition}[theorem]{Proposition}

\newtheorem{Theorem}{Theorem}[section]

\newtheorem{Corollary}[Theorem]{Corollary}

\newcommand{\R}{\mathbb R}
\newcommand{\C}{\mathbb C}
\newcommand{\N}{\mathbb N}
\newcommand{\X}{\mathbb X}
\newcommand{\Y}{\mathbb Y}
\newcommand{\Xs}{\mathscr X}
\newcommand{\loc}{\mathrm{loc}}
\newcommand{\dd}{\,\mathrm dx}

\newcommand{\ds}{\,\mathrm ds}
\newcommand{\diver}{\operatorname{div}}
\newcommand{\curl}{\operatorname{curl}}
\newcommand{\BC}{\operatorname{BC}}

\newcommand{\Dom}{\mathcal D}

\newcommand{\Rea}{\operatorname{Re}}
\newcommand{\diag}{\operatorname{diag}}

\newcommand{\ip}[2]{\left\langle #1,#2\right\rangle}
\newcommand{\norm}[1]{\left\lVert #1\right\rVert}

\newcommand{\Ub}{\overline U}
\newcommand{\BS}{\mathcal V}
\newcommand{\ca}{\mathcal A}
\newcommand{\cb}{\mathcal B}
\newcommand{\cp}{\mathcal P}

\newcommand{\crr}{\mathcal R}

\title[Stability and optimal decay of asymmetric Burgers vortices]{Three-Dimensional Stability of Weakly Asymmetric Burgers Vortices: A Circulation-Independent Asymmetry Threshold and Optimal Decay}

\author{Jingyu Jia}
\address{School of Mathematical Sciences, Shandong University, Jinan, Shandong 250100, China}
\email{202311867@mail.sdu.edu.cn}

\author{Tao Tao}
\address{School of Mathematics, Shandong University, Jinan, Shandong 250100, China}
\email{taotao@sdu.edu.cn}

\begin{document}
\begin{abstract}
	Burgers vortices are stationary solutions of the Navier--Stokes equations in a background straining flow and serve as a model for vortex tubes in three-dimensional turbulence. In this paper, we consider weakly asymmetric strains and establish the local nonlinear asymptotic stability of the corresponding Burgers vortices against small three-dimensional perturbations for arbitrary circulation. Our proof exploits the specific axial structure of the linearized operator, combining a two-dimensional resolvent analysis with decay estimates for axial derivatives to control the full three-dimensional evolution. A first-moment decomposition isolates the slow modes, and Gaussian resolvent bounds are transferred to polynomially weighted spaces. This yields nonlinear stability with the endpoint exponential decay rate $(1-\lambda)/2$, where $\lambda$ denotes the strain asymmetry, without any additional polynomial factor in time. An exact nonlinear first-moment law demonstrates that this rate is optimal over the full admissible perturbation class, including at zero circulation.
\end{abstract}
\maketitle

\pagenumbering{roman}
\tableofcontents
\clearpage
\pagenumbering{arabic}
\section{Introduction}\label{sec:introduction}
Turbulence is one of the most enduring unsolved problems in classical physics. Within the chaotic, multi-scale dynamics of turbulent flows, coherent structures emerge and play a fundamental role in energy transfer and dissipation. Among the most ubiquitous are vortex tubes, elongated regions of concentrated vorticity. Often described as the ``sinews of turbulence,’’ they are thought to be a primary channel through which energy cascades from large scales to the smallest scales, where viscous effects dominate. Understanding the mechanisms that govern their formation, persistence, and stability is therefore central to any comprehensive theory of turbulence.

Burgers vortices are a basic model for the persistence of concentrated vorticity in a straining flow. In the classical solution of Burgers~\cite{burgers1948}, vorticity amplification by axial stretching balances viscous diffusion, producing a stationary vortex with a Gaussian cross-section. Unequal transverse compression rates lead to asymmetric Burgers vortices, whose profiles are generally not explicit. We denote the strain asymmetry by $\lambda\in[0,1)$ and the circulation by $\alpha\in\mathbb R$. How these parameters affect existence, stability, and decay has motivated complementary asymptotic, numerical, and mathematical approaches, surveyed systematically by Gallay and Maekawa~\cite{GMreview}.

Early work showed that the interaction between strain and rotation is more subtle than a simple deformation of the circular core. Robinson and Saffman~\cite{robinson1984} combined perturbative calculations with numerical constructions of nonsymmetric vortices and studied the two-dimensional linear stability of the symmetric vortex. Moffatt, Kida, and Ohkitani~\cite{moffatt1994} developed a large-Reynolds-number expansion for a vortex in non-axisymmetric strain, identifying the leading axisymmetric structure and the correction induced by the asymmetry. Prochazka and Pullin then studied the axisymmetric spectrum~\cite{ProchazkaPullin1995} and the structure and stability of nonsymmetric vortices~\cite{ProchazkaPullin1998}. Their computations indicated two-dimensional stability over broad parameter ranges and showed how the vortex profile changes with strain and circulation.

The rigorous two-dimensional theory grew out of the study of the vorticity equation in self-similar variables. Giga and Kambe~\cite{GK1988} proved convergence to a Gaussian vortex under smallness assumptions and linked this behavior to the formation of Burgers vortices. Gallay and Wayne~\cite{GWglobal} later proved convergence to the Lamb--Oseen vortex for arbitrary integrable initial vorticity and arbitrary circulation; in rescaled variables, this is global attraction to the axisymmetric Burgers profile within the two-dimensional class. Their earlier invariant-manifold analysis~\cite{gallay2002} identified finite-dimensional structures governing higher-order asymptotics in weighted spaces. Together, these results explain why spatial moments and the associated low modes are central to sharp decay estimates.

For asymmetric vortices, existence and stability require a further analysis of the stationary profile. Gallay and Wayne~\cite{GW} constructed weakly asymmetric Burgers vortices for every circulation, with an asymmetry range independent of $\alpha$, justified their large-circulation asymptotics, and established local stability under two-dimensional perturbations. Maekawa~\cite{MaekawaHighRe} extended the large-circulation construction to asymmetries below $1/2$, and in a subsequent work~\cite{MaekawaAll} proved existence for every $\lambda\in[0,1)$ and every $\alpha\in\mathbb R$, together with large-circulation expansions in suitable weighted spaces. This general existence theorem does not by itself yield a stability theorem on the same parameter range.

The effect of rapid rotation has also been quantified at the level of linear and nonlinear decay. Maekawa~\cite{MaekawaSpectral2011} showed that, on an appropriate invariant subspace, the spectrum of the linearization at the axisymmetric Burgers vortex moves farther into the stable half-plane as the circulation increases. For the related Lamb--Oseen problem in self-similar variables, Li, Wei, and Zhang~\cite{LiWeiZhang2020} established sharp spectral and pseudospectral bounds, and Gallay~\cite{GallayEnhanced2018} used resolvent estimates to prove nonlinear enhanced dissipation and axisymmetrization. These works separate the slow modes associated with symmetries from components that decay faster under rapid rotation, and they show the importance of resolvent bounds for non-self-adjoint vortex operators.

Three-dimensional stability adds a further difficulty, because perturbations may vary along the vortex axis. Rossi and Le Diz`es~\cite{rossi1997} showed that the discrete temporal spectrum of a stretched vortex is associated with axially independent perturbations, so point-spectrum information alone does not describe the evolution of general three-dimensional disturbances. Schmid and Rossi~\cite{schmid2004} studied this evolution using time-dependent transformations and direct--adjoint computations. They found asymptotic damping together with potentially large transient amplification, which shows why long-time stability and short-time growth must be treated separately.

The rigorous three-dimensional theory first covered the small-circulation regime. Gallay and Wayne~\cite{GWlow} proved nonlinear asymptotic stability, allowing a change of the limiting circulation, for both symmetric and asymmetric Burgers vortices at sufficiently small Reynolds number. Gallay and Maekawa~\cite{GM3D} later established local nonlinear stability of the axisymmetric Burgers vortex for arbitrary circulation, using decay of axial derivatives to connect the three-dimensional evolution to a vectorial two-dimensional problem. They also observed that a perturbation argument gives three-dimensional stability of weakly asymmetric vortices at each fixed circulation, provided the permitted asymmetry is chosen small enough depending on that circulation. Stability of the axisymmetric vortex at arbitrary circulation and a common asymmetry threshold for the non-axisymmetric family are therefore distinct issues.

Recent work has also explored larger asymmetries at large circulation. In a preprint, Salkimbayev~\cite{Salkimbayev2026} reports a computer-assisted linear stability result beyond the half-strain threshold. This addresses a complementary parameter regime and a weighted linear problem. 
We prove three-dimensional nonlinear asymptotic stability of
weakly asymmetric Burgers vortices for arbitrary circulation.
A single asymmetry threshold applies to every circulation and
every finite polynomial weight exponent $m>3$. For small
perturbations that preserve the vortex circulation, the solution
decays at the endpoint exponential rate $(1-\lambda)/2$ in the
global perturbation norm, without a polynomial factor in time.
An exact nonlinear first-moment law shows that this rate is
optimal over the full admissible class, including at zero
circulation.

We work in nondimensional variables, with the kinematic
viscosity normalized to one. We begin by introducing the notation required to state the main result, starting with the three-dimensional incompressible Navier-Stokes equations:
\begin{equation}
	\partial_t V + (V \cdot \nabla)V = \Delta V - \nabla P, \quad \nabla \cdot V =0,
	\label{eq: ns equation}
\end{equation}
where \(V = V(x, t) \in \mathbb{R}^3\) denotes the velocity field, \(P = P(x, t) \in \mathbb{R}\) is the pressure field, and \(x = (x_1, x_2, x_3)^T \in \mathbb{R}^3\) is the spatial variable. To model tubular vortices, we decompose the velocity field as follows:
\begin{equation}
	V(x, t) = V^s(x) + U(x,t),
	\label{v decomposition}
\end{equation} 
where \(V^s(x)\) is a stationary straining flow of the form
\begin{equation}\label{vs expression}
	V^s(x)
	=
	\begin{pmatrix}
		-a_1x_1\\
		-a_2x_2\\
		x_3
	\end{pmatrix}
	=Mx,
	\qquad
	M=\operatorname{diag}(-a_1,-a_2,1),
\end{equation}
where
\[
a_1=\frac{1+\lambda}{2},
\qquad
a_2=\frac{1-\lambda}{2},
\qquad
A=\operatorname{diag}(a_1,a_2).
\]
For further convenience, instead of considering the evolution of \(V\), we prefer working with the vorticity field \(\omega = \nabla \times U\). Taking the curl of \eqref{eq: ns equation} and using \eqref{v decomposition} and \eqref{vs expression}, we obtain the evolution equation for \(\omega\):
\begin{equation}
	\partial_t \omega + (U, \nabla)\omega - (\omega, \nabla)U = L\omega, \qquad \nabla \cdot \omega = 0, 
	\label{eq: vorticity}
\end{equation}
where \(L\) is the differential operator defined by
\begin{equation}
	L\omega = \Delta\omega - (Mx, \nabla)\omega + M\omega. 
	\label{L operator}
\end{equation}

When $\lambda=0$, the background strain is axisymmetric,
with $a_1=a_2=1/2$. In this case, the vorticity equation
\eqref{eq: vorticity} admits the classical family of stationary
Burgers vortices
\begin{equation}\label{Burgers vortex}
	\Omega_{0,\alpha}(x)=\alpha G(x_h)e_3,
	\qquad
	G(x_h)=\frac{1}{4\pi}\exp\!\left(-\frac{|x_h|^2}{4}\right),
	\qquad \alpha\in\mathbb R,
\end{equation}
where $x_h=(x_1,x_2)$ and $e_3=(0,0,1)^{\mathsf T}$.
The associated rotational velocity is given by the
Biot--Savart law:
\begin{equation}\label{burgers velocity}
	U_{0,\alpha}(x_h)
	=
	\frac{\alpha}{2\pi}
	\frac{1-e^{-|x_h|^2/4}}{|x_h|^2}
	\begin{pmatrix}
		-x_2\\
		x_1\\
		0
	\end{pmatrix}.
\end{equation}

Since $\int_{\mathbb R^2}G(x_h)\,\mathrm{d}x_h=1$, the
circulation through each horizontal cross-section is
\[
\int_{\mathbb R^2}
\bigl(\Omega_{0,\alpha}\bigr)_3(x_h,x_3)\,\mathrm{d}x_h
=\alpha.
\]
In these nondimensional variables, the circulation Reynolds
number is $|\alpha|$.

The aim of this paper is to study the counterpart of Burgers vortices when the straining flow is non-axisymmetric.
Let $\omega_{\lambda,\alpha}$ denote the weakly asymmetric
stationary branch with circulation $\alpha$ constructed
in~\cite{GW}, and write
\[
\overline\Omega(x)
=\omega_{\lambda,\alpha}(x_h)e_3,
\qquad
\overline U(x)
=\bigl(u_{\lambda,\alpha}(x_h),0\bigr),
\]
where $u_{\lambda,\alpha}$ is the velocity obtained from
$\omega_{\lambda,\alpha}$ through the two-dimensional
Biot--Savart law.
At $\lambda=0$, the stationary branch reduces to the
axisymmetric family \eqref{Burgers vortex}--\eqref{burgers velocity},
with $(\overline\Omega,\overline U)
=(\Omega_{0,\alpha},U_{0,\alpha})$.

We consider solutions of \eqref{eq: vorticity} of the form \(\omega = \overline{\Omega} + \eta\) and \(U = \overline{U} + u\) with \(u=BS_3[\eta]\), where
$BS_3$ is the three-dimensional Biot--Savart operator.
For divergence-free perturbations, the evolution can be written as
\begin{equation}\label{eq:intro-perturbation}
	\begin{aligned}
		\partial_t\eta
		={}&\Delta\eta-(Mx)\cdot\nabla\eta+M\eta\\
		&+\nabla\times
		\bigl(
		\overline U\times\eta
		+u\times\overline\Omega
		+u\times\eta
		\bigr),~~
		\eta(0)=\eta_0.
	\end{aligned}
\end{equation}

An important issue is now to fix an appropriate function space for the admissible
perturbations. Since the Burgers vortex itself is essentially a two-dimensional flow, it is
natural to choose a functional setting which allows for perturbations in the same class,
but we also want to consider more general ones.

To specify the horizontal decay of the admissible perturbations, we first introduce two-dimensional spaces. Given \(m \in [0, \infty]\), we denote by \(\rho_m : [0, \infty) \to [1, \infty)\) the weight function defined by
\begin{equation}
	\rho_m(r) = 
	\begin{cases} 
		1 & \text{if } m = 0, \\ 
		(1 + \frac{r}{4m})^m & \text{if } 0 < m < \infty, \\ 
		e^{r/4} & \text{if } m = \infty. 
	\end{cases}
	\label{p weighted function}
\end{equation}

We introduce the weighted \(L^2\) space
\begin{equation}
	L^2(m) = \left\{ f \in L^2(\mathbb{R}^2) \ \middle|\ \|f\|_{L^2(m)}^2 = \int_{\mathbb{R}^2} |f(x_h)|^2 \rho_m(|x_h|^2) \, \mathrm{d}x_h < \infty \right\},
	\label{p weighted space}
\end{equation}
which is a Hilbert space with a natural inner product. Using H\"older's inequality, it is easy to verify that \(L^2(m) \hookrightarrow L^1(\mathbb{R}^2)\) if \(m > 1\). In that case, we also define the closed subspace
\begin{equation}
	L^2_0(m) = \left\{ f \in L^2(m) \ \middle|\ \int_{\mathbb{R}^2} f(x_h) \, \mathrm{d}x_h = 0 \right\}. 
	\label{zero mass p weight space}
\end{equation}

Next, we define the three-dimensional space \(X(m)\) as the set of all \(\phi : \mathbb{R}^3 \to \mathbb{R}\) for which the map \(x_h \mapsto \phi(x_h, x_3)\) belongs to \(L^2(m)\) for any \(x_3 \in \mathbb{R}\), and is a bounded and continuous function of \(x_3\). In other words, we set
\begin{equation}
	X(m) = BC(\mathbb{R}; L^2(m)), \qquad X_0(m) = BC(\mathbb{R}; L^2_0(m)), 
	\label{time-space space}
\end{equation}
where “\(BC(\mathbb{R}; Y)\)” denotes the space of all bounded and continuous functions from \(\mathbb{R}\) into \(Y\). Both \(X(m)\) and \(X_0(m)\) are Banach spaces equipped with the norm
\begin{equation}
	\|\phi\|_{X(m)} = \sup_{x_3 \in \mathbb{R}} \|\phi(\cdot, x_3)\|_{L^2(m)}. 
	\label{norm of X(m)}
\end{equation}
We study the stability of a fixed Burgers vortex with circulation
$\alpha$. The perturbations are taken in the space
\begin{equation}\label{product space}
	\mathbb X(m)=X(m)\times X(m)\times X_0(m).
\end{equation}
The condition $\eta_3\in X_0(m)$ means that the perturbation
has zero circulation on every horizontal slice. Hence the
total vorticity $\omega=\overline\Omega+\eta$ has circulation
$\alpha$ on each such slice.
We impose the incompressibility constraint by setting
\[
\mathbb Y(m)
=\left\{
\eta\in\mathbb X(m):
\nabla\cdot\eta=0
\text{ in }\mathcal D'(\mathbb R^3)
\right\}.
\]

Following the axial framework of~\cite{GM3D}, we give the definition of local topology.
\begin{definition}[Initial-value topology]\label{def:topology}
	The locally convex topology \(\tau_{\mathrm{loc}}\) on \(\mathbb{X}\) is generated by
	\[
	p_R(f) = \sum_{j=1}^3 \sup_{|z| \le R} \|f_j(\cdot, z)\|_E, \qquad R > 0.
	\]
\end{definition}
We denote by $\mathbb X_{\mathrm{loc}}(m)$ the space $\mathbb X(m)$ equipped with the topology $\tau_{\mathrm{loc}}$. In general, the axial Ornstein--Uhlenbeck evolution is not strongly continuous in the global norm of $\mathrm{BC}(\mathbb R;E)$, so we formulate continuity at the initial time in the topology $\tau_{\mathrm{loc}}$. All stability estimates, however, are stated in the global norm of $\mathbb X(m)$.

We are now able to formulate our main result:

\begin{theorem}[Uniform small-asymmetry threshold and optimal decay]
	\label{thm:main}
	There exists $\lambda_*\in(0,1)$ such that for every finite
	$m>3$, every $\lambda\in[0,\lambda_*]$, and every
	$\alpha\in\mathbb R$, there exist constants
	\[
	C=C(\lambda,\alpha,m)\ge1,
	\qquad
	\delta=\delta(\lambda,\alpha,m)>0
	\]
	with the following property.
	
	For every $\eta_0\in\mathbb Y(m)$ satisfying
	\[
	\|\eta_0\|_{\mathbb X(m)}\le\delta,
	\]
	equation~\eqref{eq:intro-perturbation} has a unique global mild
	solution in the locally norm-bounded, locally axially continuous
	solution class described above.
	Moreover, $\eta(t)\in\mathbb Y(m)$ for all $t\ge0$,
	\[
	\eta\in
	L^\infty(\mathbb R_+;\mathbb X(m))
	\cap C([0,\infty);\mathbb X_{\mathrm{loc}}(m)),
	\]
	and
	\begin{equation}\label{eq:mainbound}
		\|\eta(t)\|_{\mathbb X(m)}
		\le C e^{-\nu t}\|\eta_0\|_{\mathbb X(m)},
		\qquad
		t\ge0,
		\qquad
		\nu=\frac{1-\lambda}{2}.
	\end{equation}
	
	The exponent $\nu$ is optimal for the full class of admissible
	small perturbations.
	More precisely, for each fixed admissible triple
	$(\lambda,\alpha,m)$ and every $\mu>\nu$, no estimate of the form
	\[
	\|\eta(t)\|_{\mathbb X(m)}
	\le \widetilde C e^{-\mu t}
	\|\eta_0\|_{\mathbb X(m)},
	\qquad t\ge0,
	\]
	can hold for all sufficiently small initial data in
	$\mathbb Y(m)$, even if $\widetilde C$ and the smallness radius
	are allowed to depend on $(\lambda,\alpha,m,\mu)$.
	This optimality statement includes $\alpha=0$.
\end{theorem}
\begin{remark}
	The order of the quantifiers in Theorem~\ref{thm:main} is crucial: a single asymmetry threshold \(\lambda_*\) is valid for all circulations and all finite weights \(m>3\). 
	While the stability constant \(C\) and the perturbation radius \(\delta\) may depend on the fixed parameters, the theorem establishes local nonlinear stability at every circulation without requiring a circulation-independent basin of attraction. 
	Furthermore, Estimate~\eqref{eq:mainbound} simultaneously implies Lyapunov stability and exponential attraction in the global \(\mathbb X(m)\) norm.
\end{remark}
\begin{remark}
The optimal decay rate \(\nu=(1-\lambda)/2\) is achieved at the endpoint without an additional polynomial factor in time, and its value is dictated by the weaker transverse compression. 
To illustrate the obstruction to faster decay, consider the 
invariant class of axially independent vertical perturbations
\[
\eta(x_h,x_3,t)=(0,0,f(x_h,t)),
\qquad
\int_{\mathbb R^2}f(x_h,t)\,\mathrm dx_h=0.
\]
For these perturbations, the nonlinear evolution obeys the 
exact first-moment identities
\[
\mathcal M_j(f(t))
=e^{-a_jt}\mathcal M_j(f(0)),
\qquad
\mathcal M_j(f)
=\int_{\mathbb R^2}x_jf(x_h)\,\mathrm dx_h,
\qquad j=1,2.
\]
The first moments are continuous linear functionals on the 
weighted space under consideration. 
Taking
\[
f(0)=-\varepsilon\partial_2G,
\qquad \varepsilon\ne0,
\]
yields arbitrarily small admissible data with 
\(\mathcal M_2(f(0))=\varepsilon\).
Consequently,
\[
|\varepsilon|e^{-\nu t}
=|\mathcal M_2(f(t))|
\le C_m\|f(t)\|_{L^2(m)},
\]
which rules out every larger exponent for the full perturbation 
class. 
This argument applies at zero circulation as well. 
While perturbations satisfying additional moment cancellation 
conditions may exhibit faster decay, the optimality assertion 
concerns the complete class in Theorem~\ref{thm:main}.
\end{remark}

We now outline the proof of Theorem~\ref{thm:main} and explain the main
ideas. As in the two-dimensional theory, the essential task is to establish
quantitative stability estimates for the linearization around the stationary
vortex. The new difficulties are the lack of radial symmetry, the presence of
arbitrarily large circulation, and the coupling of horizontal and axial
perturbations. Our argument proceeds in the following steps.

\medskip
\noindent\textbf{Step 1: Stability of the horizontal vector operator.}
We use the weakly asymmetric stationary vortex from~\cite{GW}, without
imposing a smallness condition on the circulation $\alpha$. The horizontal
linearization is the operator $H_b$ introduced in
Section~\ref{sec:horizontal}. Since the stationary velocity is generally
nonradial and its size need not be small, a direct coercive estimate for
$H_b$ is insufficient. Instead, we exploit the divergence and
stream-function identities in~\eqref{eq:intertwine} to locate its point
spectrum. A compactly supported absorbing potential and a Fredholm
argument then exclude the entire half-plane $\{\Re z\ge-1/2\}$ from the
spectrum and give the uniform resolvent estimate
\[
\sup_{\Re z\ge-1/2}
\bigl\|(z-H_b)^{-1}\bigr\|_{E^2\to E^2}<\infty.
\]
The Gearhart--Pr\"uss theorem consequently yields
\[
\bigl\|e^{tH_b}\bigr\|_{E^2\to E^2}
\le C_h e^{-t/2},\qquad t\ge0.
\]
The constants may depend on $\alpha$, whereas the small-asymmetry
threshold does not.

\medskip
\noindent\textbf{Step 2: Stability of the vertical scalar operator.}
The scalar block $A_3$ is first treated in a Gaussian weighted space,
where a coercive form estimate is available on the zero-mass subspace.
To transfer this estimate to the finite polynomial weight, we introduce
a rank-one mass damping and use a triangular Gaussian extension. The
remaining perturbation is regularizing from the polynomial space into
the Gaussian space, which permits a variation-of-constants argument.
This gives the initial bound
\[
\bigl\|e^{tA_3}\bigr\|_{E_0\to E_0}
\le C_v e^{-t/4},\qquad t\ge0.
\]
Unlike an estimate obtained by simply treating the entire stationary
velocity as a small perturbation, the argument remains valid for each
$\alpha\in\mathbb R$ under a single asymmetry threshold.

\medskip
\noindent\textbf{Step 3: Three-dimensional linear stability.}
The full linearized operator is decomposed into a two-dimensional
principal part and a genuinely three-dimensional remainder:
\[
\mathcal B=\mathcal L_0+\mathcal R,
\qquad
\mathcal L_0
=\operatorname{diag}(H_b,A_3)
+(\partial_z^2-z\partial_z).
\]
The principal part inherits the preceding two-dimensional decay
estimates. The key additional observation is the commutator identity
\[
[\partial_z,\mathcal B]=-\partial_z,
\qquad
\partial_z^k\widehat S(t)f
=e^{-kt}\widehat S(t)\partial_z^k f,
\]
which expresses the damping of axial derivatives by the background
stretching. Combined with axial smoothing and the remainder bound
\eqref{eq:Rest}, this allows us to control $\mathcal R$ and prove
\[
\bigl\|\widehat S(t)\bigr\|_{\mathbb X(m)\to\mathbb X(m)}
\le M e^{-t/4},\qquad t\ge0.
\]
We first construct this evolution on the full product space
$\mathbb X(m)$ and only then restrict it to the divergence-free
subspace $\mathbb Y(m)$; the two pieces of the operator do not
separately preserve the divergence constraint.

\medskip
\noindent\textbf{Step 4: Nonlinear stability.}
To handle the quadratic vorticity nonlinearity, we construct a curl
propagation operator $Z(t)$ extending $S(t)\operatorname{curl}$.
Its essential short-time estimate is
\[
\bigl\|Z(t)F\bigr\|_{\mathbb X(m)}
\le C_Ze^{-t/4}(1-e^{-t})^{-3/4}
\|F\|_{\mathscr X_{4/3}},\qquad t>0.
\]
The power $3/4<1$ makes the singularity integrable and allows us to
solve the nonlinear Duhamel equation
\[
\eta(t)=S(t)\eta_0
+\int_0^t Z(t-s)\mathscr Q(\eta(s))\,\mathrm ds,
\qquad
\mathscr Q(\eta)=(BS_3[\eta])\times\eta.
\]
by a contraction argument in an exponentially weighted trajectory
space. Since axial translations need not be strongly continuous in
the global supremum norm, we use the local axial topology at the
initial time and justify the integral in local Banach spaces. This
establishes global existence, uniqueness, and the baseline decay rate
for sufficiently small perturbations.

\medskip
\noindent\textbf{Step 5: The optimal rate and its sharpness.}
To improve the baseline decay to the optimal exponent
$\nu=(1-\lambda)/2$, we isolate the first horizontal moments of the
scalar evolution. The exact identities
\[
p_j(A_3f)=-a_jp_j(f),\qquad
p_j(f)=\int_{\mathbb R^2}x_j f(x_h)\,\mathrm dx_h,
\qquad j=1,2,
\]
separate the slowly decaying modes from the subspace where both the
mass and first moments vanish. On the subspace of zero mass and vanishing first moments, the polynomially weighted evolution decays at a rate $\sigma>\nu$.
Transferring the resulting bound to the full three-dimensional
evolution yields
\[
\bigl\|\widehat S(t)\bigr\|_{\mathbb X(m)\to\mathbb X(m)}
\le C_\sharp e^{-\nu t},\qquad t\ge0,
\]
with no additional polynomial factor in time. The curl estimate and
the nonlinear fixed point can then be repeated with the endpoint
weight $e^{\nu t}$, proving~\eqref{eq:mainbound}. Finally, for
axially independent vertical perturbations, the nonlinear evolution
satisfies the exact first-moment law
\[
p_2(f(t))=e^{-a_2t}p_2(f(0)),
\qquad a_2=\nu.
\]
Choosing initial data with $p_2(f(0))\ne0$ rules out any faster
uniform exponential decay and proves sharpness.

\medskip
\noindent\textbf{Organization of the paper.}
Section~\ref{sec:inputs} introduces the stationary vortex,
establishes its spatial regularity, and develops the operator
realizations and weighted smoothing estimates used below.
Section~\ref{sec:linear} proves the horizontal and scalar
linear stability estimates and constructs the full
three-dimensional evolution.
Section~\ref{sec:curlprop} constructs the curl propagation
operator and establishes nonlinear existence, uniqueness,
and stability with a baseline decay rate.
Finally, Section~\ref{sec:optimal} proves the optimal decay
rate, establishes its sharpness, and records the dependence
of the constants on the parameters.

\section{The stationary profiles and weighted semigroup estimates}
\label{sec:inputs}
This preliminary section collects fundamental estimates for the main theorem \ref{thm:main}, establishing the spatial regularity of the asymmetric Burgers vortex and the Gaussian realization and exponential decay of the scalar linearization operator $A_3$.

\subsection{Regularity and decay of asymmetric vortices}
\label{sec:asymmetric Burgers vortex}
We use the weakly asymmetric stationary vortices constructed
by Gallay and Wayne in \cite[Proposition~4.2]{GW}.
More precisely, there exist constants
$\lambda_B\in(0,1)$ and $C_B>0$, independent of the
circulation $\alpha\in\mathbb R$, such that the following
statement holds for every
$0\leq\lambda\leq\lambda_B$ and every $\alpha\in\mathbb R$.

Set
\[
G(x_h)=\frac{1}{4\pi}e^{-|x_h|^2/4},
\qquad
F=L^2(\mathbb R^2,G^{-1}\,\mathrm dx_h),
\qquad
\ell(f)=\int_{\mathbb R^2}f(x_h)\,\mathrm dx_h,
\]
and define
\[
F_0=\{f\in F:\ell(f)=0\},
\qquad
F^1=\{f\in F:\nabla_h f\in F^2\},
\qquad
F^2=F\times F,
\]
with
\[
\|f\|_{F^1}
=
\|f\|_F+\|\nabla_h f\|_{F^2}.
\]
There exists a stationary vortex
$\omega=\omega_{\lambda,\alpha}$ of the form
\[
\omega=\alpha G+w,
\qquad
w\in F_0\cap F^1,
\]
satisfying
\begin{equation}\label{eq:branch}
	\|w\|_{F^1}
	\leq
	C_B\lambda\frac{|\alpha|}{1+|\alpha|}.
\end{equation}
Its velocity $u=K_2*\omega$ and vorticity satisfy
\begin{equation}\label{eq:branch-stationary}
	\Delta_h\omega
	+(Ax_h-u)\cdot\nabla_h\omega+\omega=0
	\quad\text{in }\mathcal D'(\mathbb R^2),
	\qquad
	\ell(\omega)=\alpha.
\end{equation}
The existence of this branch and estimate
\eqref{eq:branch} are taken from
\cite[Proposition~4.2]{GW}.
The additional spatial regularity needed in this paper
is established directly below from
\eqref{eq:branch-stationary} and the Gaussian weighted
bound.

\begin{lemma}[Regularity and decay of the asymmetric Burgers vortex]\label[lemma]{lem:branch-regularity}
	Fix $\alpha\in\mathbb{R}$ and $0\leq\lambda\leq\lambda_B$.
	Let $\omega=\alpha G+w$ and $u=K_2*\omega$ be the stationary
	vorticity and velocity described in \eqref{eq:branch-stationary}. Then
	\[
	\omega\in\mathcal{S}(\mathbb{R}^2),
	\qquad
	\nabla\omega\in\mathbf{F}.
	\]
	Moreover, $u$ is smooth, and
	$u$ and all its derivatives are bounded. In particular,
	\[
	\rho_m(|x|^2)^{1/2}\nabla\omega
	\in L^4(\mathbb{R}^2)^2,
	\qquad \forall\ 0\leq m<\infty.
	\]
\end{lemma}
\begin{proof}
	Based on the above $\omega$ assumption, it can be easily deduced that $\omega\in L^1\cap H^1$ and
	$\nabla\omega\in \mathbb{F}$. According to the two-dimensional Biot-Savart law, it can be known that
	\begin{equation}\label{eq:kernel-elementary}
		\norm{K_2*f}_4\le C(\norm{f}_1+\norm{f}_2),\qquad
		\norm{K_2*g}_\infty\le C(\norm{g}_1+\norm{g}_4).
	\end{equation}
	Since \(H^1(\R^2)\hookrightarrow L^4(\R^2)\) and \(f\in L^1\cap L^2\), it follows immediately that \(u=K_2*\omega\) is bounded.
	
	We next use the stationary equation to obtain higher derivatives
	of $\omega$. Write \eqref{eq:branch-stationary} as
	\[
	\Delta_h\omega=(u-Ax_h)\cdot\nabla_h\omega-\omega.
	\]
	The Gaussian weight in the initial $F^1$ bound implies that,
	for every $N\geq0$,
	\[
	\|\omega\|_{H^1(B_2(x_0))}
	\leq C_N(1+|x_0|)^{-N}.
	\]
	Here and below, the constants $C$ are independent of $x_0$.
	Since $u$ is bounded and $Ax_h$ grows at most linearly, the interior elliptic estimate
	\cite[Section~6.3.1, Theorem~1]{Evans} gives
	\[
	\|\omega\|_{H^2(B_1(x_0))}
	\leq C(1+|x_0|)
	\|\omega\|_{H^1(B_2(x_0))}
	\leq C_N(1+|x_0|)^{-N}.
	\]
	where in the last inequality, we use the preceding bound with
	exponent $N+1$ to absorb the factor $1+|x_0|$.
	
	The embedding $H^2(B_1)\hookrightarrow W^{1,4}(B_1)$ \cite[Corollary~9.14]{Brezis},
	together with a covering by unit balls with bounded overlap,
	shows that $\nabla\omega$ belongs to $L^1\cap L^4$, even after
	multiplication by any polynomial weight. Consequently,$\partial_j u=K_2*(\partial_j\omega)$
	and \eqref{eq:kernel-elementary} imply that the first
	derivatives of $u$ are bounded.
	
	We can now repeat the same argument. Suppose that for some integer \(k\ge 0\), we have $\sum_{|\beta|\leq k}
	\|\partial^\beta u\|_{L^\infty(\mathbb{R}^2)}
	<\infty$, and that
	\[
	\forall N\in\mathbb{N}_0,\quad
	\exists C_{k,N}<\infty,\quad
	\|\omega\|_{H^{k+1}(B_1(x_0))}
	\leq C_{k,N}(1+|x_0|)^{-N}
	\quad\text{for all }x_0\in\mathbb{R}^2.
	\] Differentiating the stationary equation \eqref{eq:branch-stationary} up to
	$k$ times and applying the interior Laplacian estimate yields
	\[
	\|\omega\|_{H^{k+2}(B_1(x_0))}
	\leq C_k(1+|x_0|)
	\|\omega\|_{H^{k+1}(B_2(x_0))}.
	\]
	The decay bounds on $B_2(x_0)$ follow by covering it with
	finitely many unit balls, and the additional factor $1+|x_0|$
	can again be absorbed by using a higher decay exponent.
	Thus the local $H^{k+2}$ norms have the same rapid decay.
	
	Applying the two-dimensional embedding $H^{k+2}(B_1)\hookrightarrow W^{k+1,4}(B_1)$, and \eqref{eq:kernel-elementary} to these derivatives
	shows that the derivatives of $u$ through order $k+1$ are
	bounded. This completes the induction.
	
	Finally, local Sobolev embedding at arbitrarily high orders
	shows that every derivative of $\omega$ has rapid pointwise
	decay. Hence $\omega\in\mathcal{S}(\mathbb{R}^2)$.
	Since $\rho_m(|x|^2)^{1/2}$ grows at most polynomially for
	finite $m$, the stated weighted $L^4$ estimate follows.
\end{proof}

We fix
\begin{equation}\label{eq:lambdastar}
	\lambda_*:=\min\{\lambda_B,\lambda_G,1/4\},\qquad \gamma:=1/4.
\end{equation}
Here $\lambda_G>0$ is the circulation-independent constant
defined in \eqref{eq:A3-Gaussian-threshold} in the proof of
Proposition~\ref{prop:scalar-resolvent}.

The symbol $\gamma=1/4$ denotes the baseline rate throughout; the
optimal rate $\nu=(1-\lambda)/2$ is introduced separately in
Section~\ref{sec:optimal}. Every other constant below may depend on the fixed triple
$(\lambda,\alpha,m)$. The two numerical margins that are independent of
$\alpha$ are
\begin{equation}\label{eq:margins}
	ma_2-\frac12>\frac58,\qquad
	\frac12+(m+1)a_2>2.
\end{equation}

Since $\omega$ is Schwartz, its Biot--Savart velocity obeys
\begin{equation}\label{eq:udecay}
	u,\nabla u\in L^\infty,\qquad
	u(x_h)=O(|x_h|^{-1}),\quad \nabla u(x_h)=O(|x_h|^{-2}).
\end{equation}
For completeness, we split the convolution at $|y|=|x|/2$. In the inner region,
we expand the kernel and use the integrable moments of $\omega$. In the remaining
region, we use its rapid decay, treating the local singularity by subtracting a
constant when differentiating the kernel. This gives \eqref{eq:udecay}, with
constants depending on the fixed profile.

For finite $m>3$, 
set
\[
E=L^2(\mathbb R^2,\rho_m \,dx_h),
\qquad
\|f\|_E^2
=
\int_{\mathbb R^2}|f(x_h)|^2\rho_m(x_h)\,dx_h,
\]
and
\[
E_0=\{f\in E:\ell(f)=0\}.
\]
The weight satisfies
\begin{equation}\label{eq:weights}
	\begin{gathered}
		|\nabla_h\log \rho_m(x_h)|
		\leq\frac{C_m}{1+|x_h|},
		\qquad
		\frac{\Delta_h\rho_m(x_h)}{\rho_m(x_h)}
		=
		O((1+|x_h|)^{-2}),\\
		\rho_m(x_h+y_h)^{1/2}
		\leq
		C_m(1+|y_h|)^m\rho_m(x_h)^{1/2}.
	\end{gathered}
\end{equation}
By Cauchy--Schwarz,
\begin{equation}\label{eq:embeddings}
	\|f\|_1+\|f\|_2\leq C_m\|f\|_E,
	\qquad
	|\ell(f)|
	\leq
	\left(\frac{4\pi m}{m-1}\right)^{1/2}\|f\|_E.
\end{equation}
Since $m>2$,
\begin{equation}\label{eq:moment}
	\|(1+|x_h|)f\|_1\leq C_m\|f\|_E.
\end{equation}
Moreover, $F\hookrightarrow E$ continuously
for every finite $m$.

\subsection{Maximal realizations and weighted smoothing estimates}
\label{sec:domains}
This section establishes the maximal distribution representation of local weighted operators. We prove that it generates a \(C_0\)-semigroup, establish gradient graph norm estimates, and derive short-time weighted smoothing bounds.
\begin{lemma}[Maximal realization and semigroup generation]\label[lemma]{lem:domain}
Consider the differential expression
\[
\mathcal{L}f
=\Delta_h f+\beta(x)\cdot\nabla_h f+C(x)f,
\qquad
\beta(x)=B_0x_h+c(x),
\]
where
\[
B_0\in\mathbb{R}^{2\times2},
\qquad
c\in C^\infty(\mathbb{R}^2;\mathbb{R}^2),
\qquad
C\in C^\infty(\mathbb{R}^2;\mathbb{C}^{N\times N}),
\]
and $c$, $\nabla c$, and $C$ are bounded.

Define the maximal distribution realization by
\[
\begin{aligned}
	D(\mathcal{L}_{\max})
	=
	\left\{
	f\in E^N:
	\mathcal{L}f\in E^N
	\text{ in }\mathcal{D}'(\mathbb{R}^2)^N
	\right\},\quad \mathcal{L}_{\max}f=\mathcal{L}f.
\end{aligned}
\]
Then
\[
\mathcal{L}_{\max}
=
\overline{
	\mathcal{L}\big|_{C_c^\infty(\mathbb{R}^2)^N}
},
\]
and $\mathcal{L}_{\max}$ generates a $C_0$-semigroup $U(t)$
on $E^N$. Moreover,
\begin{equation}\label{eq:gradgraph}
 \norm{\nabla f}_E^2\le C\norm{f}_E^2+
 C\norm{\mathscr Lf}_E\norm{f}_E,
 \qquad f\in\Dom_{\max}(\mathscr L).
\end{equation}
Define the Hermitian matrix-valued energy potential
\[
V_\rho(x)
=
\frac{C(x)+C(x)^*}{2}
+
\frac12
\left(
\frac{\Delta \rho_m(x)}{\rho_m(x)}
-\operatorname{div}\beta(x)
-\beta(x)\cdot\nabla\log \rho_m(x)
\right)I_N.
\]
If $V_\rho(x)\leq c_0I_N$ for all $x$, then
\[
\|U(t)\|_{E^N\to E^N}\leq e^{c_0t},
\qquad t\geq0.
\]
In particular, if $V_\rho(x)\leq-\kappa I_N$ for some
$\kappa>0$, then
\[
\|U(t)\|_{E^N\to E^N}\leq e^{-\kappa t},
\qquad t\geq0.
\]
The matrix inequalities are understood in the sense of
Hermitian quadratic forms.
\end{lemma}
\begin{remark}
    In the applications below, we use $N=1$ for the scalar
operators ($P_b, Q$) and $N=2$ for the horizontal vector operators ($H_b, Q_h$).
\end{remark}
\begin{proof}
We first observe that the maximal realization is closed. Indeed, onvergence of both
$f_n$ and $\mathcal L f_n$ in $E^N$ allows passage to the
limit in distributions. It is densely defined because
its domain contains $C_c^\infty(\mathbb R^2)^N$.

Let $f\in D(\mathcal L_{\max})$. The local elliptic regularity
yields $f\in H^2_{\mathrm{loc}}(\mathbb R^2)^N$.
Choose $\chi\in C_c^\infty(\mathbb R^2)$ with
$0\le\chi\le1$, $\chi=1$ on $B_1$, and
$\operatorname{supp}\chi\subset B_2$, and set
$\chi_R(x)=\chi(x/R)$ for $R\ge1$.
Testing $\mathcal L f$ against $\chi_R^2 \rho_m \overline f$,
integrating by parts, and using Young's inequality gives
\[
  \frac12
  \int_{\mathbb R^2}\chi_R^2\rho_m|\nabla_h f|^2\,dx
  \le
  K\|f\|_E^2+\|\mathcal L f\|_E\|f\|_E,
\]
with $K$ independent of $R$. 
The Fatou's lemma now yields $\nabla_h f\in E^{2N}$ and \eqref{eq:gradgraph}.

A direct computation yields the commutator identity
\begin{equation}\label{eq:commcutoff}
	[\mathscr L,\chi_R]f
	=2\nabla_h\chi_R\cdot\nabla_h f+
	(\Delta_h\chi_R+\beta\cdot\nabla_h\chi_R)f\longrightarrow0\quad\text{in }E^N.
\end{equation}
Indeed, the coefficients of this commutator are uniformly bounded and supported on the annulus $\{R\le|x|\le2R\}$, while $f$ and $\nabla_h f$ possess integrable weighted square tails. Consequently, $\chi_R f\to f$ in the graph norm of the maximal realization $\mathcal L_{\max}$. For each fixed $R$, mollifying $\chi_R f$ yields smooth approximations in $H^2$ supported in a common compact set. Since the coefficients and the weight are bounded on this set, the approximation is also valid in the graph norm. Combining this with the closedness of $\mathcal L_{\max}$, we obtain the asserted core property.

For $f\in C_c^\infty(\mathbb R^2)^N$, integration by parts
gives
\begin{equation}\label{eq:localenergy}
 \Rea\ip{\mathscr Lf}{f}_{E^d}
 =-\norm{\nabla_h f}_E^2+\int \ip{V_\rho(x)f(x)}{f(x)}_{\C^d}\rho_m(x)\dd,
\end{equation}
The graph estimate implies convergence of gradients along
any approximation in the graph norm. Since $V_\rho$ is bounded,
the core property therefore extends
\eqref{eq:localenergy} to every
$f\in D(\mathcal L_{\max})$. In particular,
\[
  \operatorname{Re}
  \langle\mathcal L_{\max}f,f\rangle_{E^N}
  \le
  -\|\nabla f\|_E^2+c_0\|f\|_E^2.
\]

It remains to prove the range condition. Fix a real $s>c_0$
and $g\in E^N$. On $H_0^1(B_R)^N$, consider the
sesquilinear form
\[
  \begin{aligned}
    a_R(v,h)
    =\int_{B_R}\rho_m\Bigg[
      \sum_{j=1}^2
       \partial_jv\cdot\overline{\partial_jh} 
      +
      \Big(
        (\nabla_h\log \rho_m-\beta)\cdot\nabla_h v
        +(sI_N-C)v
      \Big)\cdot\overline h
    \Bigg]\,dx.
  \end{aligned}
\]
For each fixed $R$, this form is continuous, and integration
by parts gives
\[
  \operatorname{Re}a_R(v,v)
  \ge
  \|\nabla_h v\|_{L^2(B_R,W)}^2
  +(s-c_0)\|v\|_{L^2(B_R,W)}^2.
\]
The Lax--Milgram theorem \cite[Corollary~5.8 and Proposition~11.28]{Brezis} therefore provides
$f_R\in H_0^1(B_R)^N$ satisfying
\[
  a_R(f_R,h)
  =
  \int_{B_R}\rho_m g\cdot\overline h\,dx
  \qquad
  \text{for all }h\in H_0^1(B_R)^N.
\]
Equivalently, $(s-\mathcal L)f_R=g$ distributionally
in $B_R$. Extending $f_R$ by zero outside $B_R$, we obtain
\[
  \|\nabla_h f_R\|_E^2
  +(s-c_0)\|f_R\|_E^2
  \le
  \|g\|_E\|f_R\|_E.
\]
These extensions are uniformly bounded in
$H^1(\mathbb R^2,W\,dx)^N$. Along a sequence $R\to\infty$,
they converge weakly to some $f$ in this space.
Testing against functions supported in a fixed compact set
and passing to the limit shows that
\[
  (s-\mathcal L)f=g
  \quad\text{in }\mathcal D'(\mathbb R^2)^N.
\]
Hence $\mathcal L f=sf-g\in E^N$, so
$f\in D(\mathcal L_{\max})$ and
\[
  \operatorname{Ran}(s-\mathcal L_{\max})=E^N.
\]

Consequently, $\mathcal L_{\max}-c_0I$ is maximal
dissipative. The Lumer--Phillips theorem \cite{LP}, in the formulation of
\cite[Theorem~II.3.15]{EN}, yields a
contraction semigroup generated by this shifted operator,
and thus
\[
  \|U(t)\|_{E^N\to E^N}
  \le e^{c_0t},
  \qquad t\ge0.
\]
Taking $c_0=-\kappa$ gives the final assertion.
\end{proof}

\begin{lemma}[Short-time weighted smoothing]\label[lemma]{lem:localsmooth}
Under the hypotheses and notation of \cref{lem:domain}, assume that
\[
  B_0=\operatorname{diag}(b_1,b_2),
  \qquad b_1,b_2\in\mathbb R,
\]
and let $U(t)$ be the $C_0$-semigroup generated by
$\mathcal L_{\max}$ on $E^N$. Then there exists $K>0$ such that
\[
  \|U(t)f\|_E
  +\sqrt{t}\,\|\nabla_h U(t)f\|_E
  \le K\|f\|_E,
  \qquad
  f\in E^N,\quad 0<t\le1.
\]
If, in addition, there exists $\kappa>0$ such that $\norm{U(t)}\le e^{-\kappa t}$, then
\begin{equation}\label{eq:U14}
 \norm{U(t)f}_1+\norm{U(t)f}_4
 \le C(1+t^{-1/4})e^{-\kappa t}\norm{f}_E,
 \qquad t>0.
\end{equation}
\end{lemma}
\begin{proof}
We decompose
\[
  \mathcal L=\mathcal L_0+\mathcal B,
  \qquad
  \mathcal L_0=\Delta+B_0x\cdot\nabla,
  \qquad
  \mathcal Bf=c(x)\cdot\nabla f+C(x)f.
\]
Since $c$ and $C$ are bounded,
\[
  \|\mathcal Bf\|_E
  \le K\bigl(\|f\|_E+\|\nabla f\|_E\bigr),
  \qquad
  f\in H^1(\mathbb R^2,\rho_m \,dx)^N.
\]
Let $\mathcal L_{0,\max}$ denote the maximal distributional
realization of $\mathcal L_0$ on $E^N$, and let $S_{\mathcal{L}_0}(t)$ be its
semigroup. The explicit Gaussian representation is
\[
  S_{\mathcal{L}_0}(t)f(x)=(k_t*f)(e^{tB_0}x),
\]
where
\[
  k_t(y)
  =
  \frac{1}{4\pi\sqrt{q_1(t)q_2(t)}}
  \exp\left(
    -\frac{y_1^2}{4q_1(t)}
    -\frac{y_2^2}{4q_2(t)}
  \right),
  \qquad
  q_j(t)=\int_0^t e^{2b_js}\,ds.
\]
For $0<t\le1$, one has $q_j(t)\asymp t$, and hence
\[
  \|(1+|\cdot|)^m k_t\|_{L^1}
  +\sqrt{t}\,
   \|(1+|\cdot|)^m\nabla k_t\|_{L^1}
  \le K.
\]
Weighted Young inequalities and the change of variables
$y=e^{tB_0}x$ therefore give
\[
  \|S_{\mathcal{L}_0}(t)f\|_E
  +\sqrt{t}\,\|\nabla S_{\mathcal{L}_0}(t)f\|_E
  \le K\|f\|_E,
  \qquad 0<t\le1.
\]

We next relate this estimate to the semigroup $U(t)$.
By \cref{lem:domain}, the domains of both $\mathcal L_{\max}$ and
$\mathcal L_{0,\max}$ are contained in
$H^1(\mathbb R^2,\rho_m\,dx)^N$. Thus $\mathcal Bf\in E^N$
whenever $f$ belongs to either domain. Since
$\mathcal Lf=\mathcal L_0f+\mathcal Bf$ distributionally,
the maximal-domain definitions imply
\[
  D(\mathcal L_{\max})=D(\mathcal L_{0,\max}).
\]
Consequently, for $f$ in this common domain, the
variation-of-constants formula gives
\[
  U(t)f
  =
  S_{\mathcal{L}_0}(t)f
  +\int_0^t S_{\mathcal{L}_0}(t-s)\mathcal B U(s)f\,ds.
\]
For $0<T\le1$, set
\[
  M_T(f)
  =
  \sup_{0<t\le T}
  \left(
    \|U(t)f\|_E
    +\sqrt{t}\,\|\nabla U(t)f\|_E
  \right).
\]
For $f\in D(\mathcal L_{\max})$, this quantity is finite
by the graph estimate in Lemma~2.3 and continuity of
$U(t)f$ in the graph norm. Moreover,
\[
  \|\mathcal B U(s)f\|_E
  \le K s^{-1/2}M_T(f),
  \qquad 0<s\le T.
\]
Applying the estimates for $S(t)$ to the
variation-of-constants formula yields
\[
  M_T(f)
  \le K\|f\|_E+K\sqrt{T}\,M_T(f),
\]
where we used
\[
  \int_0^t s^{-1/2}\,ds=2\sqrt{t},
  \qquad
  \int_0^t(t-s)^{-1/2}s^{-1/2}\,ds=\pi.
\]
Choose $T_0\in(0,1]$ sufficiently small that
$K\sqrt{T_0}\le1/2$. Absorbing the last term gives
\[
  M_{T_0}(f)\le 2K\|f\|_E.
\]
By density of $D(\mathcal L_{\max})$ in $E^N$, this estimate
extends to all $f\in E^N$. The semigroup property and
boundedness of $U(t)$ on $[0,1]$ then extend the estimate
from $(0,T_0]$ to $(0,1]$. This proves the first assertion.

Assume now that
$\|U(t)\|_{E^N\to E^N}\le e^{-\kappa t}$.
Since $m>1$ and $W\ge1$,
\[
  \|g\|_{L^1}+\|g\|_{L^2}\le K\|g\|_E.
\]
Combining these embeddings with the two-dimensional
Gagliardo--Nirenberg inequality \cite[Theorem~9.9]{Brezis}
\[
  \|g\|_{L^4}
  \le
  K\|g\|_{L^2}^{1/2}
   \|\nabla g\|_{L^2}^{1/2}
\]
and the short-time estimate gives
\[
  \|U(t)f\|_{L^1}+\|U(t)f\|_{L^4}
  \le K(1+t^{-1/4})\|f\|_E,
  \qquad 0<t\le1.
\]
Since $e^{-\kappa t}\ge e^{-\kappa}$ on this interval,
this implies \eqref{eq:U14}
for $0<t\le1$, after adjusting the constant.

For $t\ge1$, write
\[
  U(t)=U(1/2)U(t-1/2).
\]
The short-time estimate at time $1/2$ and the assumed
exponential bound imply
\[
  \begin{aligned}
    \|U(t)f\|_E+\|\nabla_h U(t)f\|_E
    &\le K\|U(t-1/2)f\|_E\\
    &\le K e^{-\kappa(t-1/2)}\|f\|_E\\
    &\le K_\kappa e^{-\kappa t}\|f\|_E.
  \end{aligned}
\]
Applying the same embeddings and interpolation inequality
proves \eqref{eq:U14} for $t\ge1$,
and completes the proof.
\end{proof}

\section{The resolvent estimates and three-dimensional linear stability}
\label{sec:linear}
This section establishes decay for the horizontal and scalar linear blocks and then for the full
three-dimensional evolution.
\subsection{The horizontal resolvent bounds and semigroup decay}
\label{sec:horizontal}
This section establishes the exponential stability of the horizontal vector block. We first derive \(L^1\) and point spectrum bounds without radial symmetry, then use resolvent estimates to prove the semigroup decays as \(e^{-t/2}\).

Set $b=-Ax_h+u$, $D_b=\Delta_h-b\cdot\nabla_h$, and
\begin{equation}\label{eq:H}
 H_bh=D_bh+(h\cdot\nabla_h)b.
\end{equation}
All operators in this section are complexified; stability for real fields
follows by restriction.

\begin{lemma}[$L^1$ eigenvalue bound]\label[lemma]{lem:L1}
Suppose that $b$ is real and smooth,
$\operatorname{div}_h b=-1$, and
$|b(x_h)|\leq C(1+|x_h|)$.
If $0\neq f\in L^1(\mathbb{R}^2)$ satisfies
\[
(D_b+c)f=zf
\quad\text{in }\mathcal{D}'(\mathbb{R}^2),
\qquad c\in\mathbb{R},
\]
then
\[
\operatorname{Re}z\leq c-1.
\]
\end{lemma}
\begin{proof}
The local elliptic regularity implies that $f$ is smooth.
The complex Kato inequality \cite{Simon} gives
\[
D_b|f|
\geq(\operatorname{Re}z-c)|f|
\quad\text{in }\mathcal{D}'(\mathbb{R}^2).
\]
Choose a nonnegative smooth cutoff
$\chi_R(x_h)=\chi(x_h/R)$, equal to one on $B_R$
and supported in $B_{2R}$.
Testing the preceding inequality against $\chi_R$ and
using $\operatorname{div}_h b=-1$, we obtain
\[
(\operatorname{Re}z-c+1)
\int_{\mathbb{R}^2}\chi_R|f|\,dx_h
\leq
\int_{\mathbb{R}^2}
|f|(\Delta_h\chi_R+b\cdot\nabla_h\chi_R)\,dx_h.
\]
The coefficient on the right is uniformly bounded
and supported in the annulus
$\{R\leq|x_h|\leq2R\}$.
Since $f\in L^1$, the right-hand side tends to zero.
Letting $R\to\infty$ gives
\[
(\operatorname{Re}z-c+1)\|f\|_1\leq0.
\]
The conclusion follows because $f\neq0$.
\end{proof}

\begin{lemma}[Point spectrum without radial symmetry]\label[lemma]{lem:point}
Let $m>2$.
Every eigenvalue of $H_b$ on $E^2$ satisfies
\[
\operatorname{Re}z\leq-1.
\]
If an eigenfunction is horizontally divergence free,
its eigenvalue satisfies
\[
\operatorname{Re}z\leq-2.
\]
\end{lemma}
\begin{proof}
The direct differentiation gives
\begin{equation}\label{eq:intertwine}
 \diver_h(H_bh)=D_b\diver_hh,\qquad
 H_b(J\nabla \psi)=J\nabla(D_b\psi-\psi).
\end{equation}
For the first identity the mixed first derivatives cancel and
$h\cdot\nabla\diver b=0$. For the second use
$(\nabla b)J+J(\nabla b)^T=(\diver b)J=-J$.

Suppose that $H_bh=zh$ with $h\neq0$.
Lemma \ref{lem:domain} gives $\nabla_hh\in E$.
The embeddings \eqref{eq:embeddings} and \eqref{eq:moment}, with $m>2$, therefore imply
\begin{equation}\label{eq:hL1}
 h,\nabla_h h\in L^1,\qquad |x_h|h\in L^1.
\end{equation}
If $d:=\operatorname{div}_h h\neq0$, then $D_bd=zd$ and \cref{lem:L1} gives
$\Rea z\le-1$.

Now suppose that $\operatorname{div}_h h=0$.
We define
\[
F_2(x_2)=\int_{\mathbb{R}}h_2(s,x_2)\,ds.
\]
By Fubini's theorem and the $W^{1,1}$ slice theorem,
$F_2\in W^{1,1}(\mathbb{R})$ and
\[
F_2'(x_2)
=
\int_{\mathbb{R}}\partial_2h_2(s,x_2)\,ds
=
-\int_{\mathbb{R}}\partial_1h_1(s,x_2)\,ds
=0.
\]
Since $F_2\in L^1(\mathbb{R})$, it follows that $F_2=0$.
Consequently,
\[
\psi(x_1,x_2)
:=
\int_{-\infty}^{x_1}h_2(s,x_2)\,ds
=
-\int_{x_1}^{\infty}h_2(s,x_2)\,ds
\]
satisfies
\begin{equation}\label{eq:streamL1}
 \norm{\psi}_1\le\int_{\R^2}|x_1||h_2(x)|\dd<\infty,
 \qquad \partial_1\psi=h_2,\quad\partial_2\psi=-h_1.
\end{equation}
Thus $h=J\nabla_h\psi$.

Using \eqref{eq:intertwine} in the eigenvalue
equation gives
\[
J\nabla_h\bigl(D_b\psi-(1+z)\psi\bigr)=0.
\]
Hence $D_b\psi-(1+z)\psi$ is constant.
This constant belongs to $L^1$, because
\[
\Delta_h\psi=\partial_1h_2-\partial_2h_1\in L^1,
\qquad
b\cdot\nabla_h\psi\in L^1,
\]
by \eqref{eq:hL1} and the linear
growth of $b$.
It must therefore be zero, so
\[
(D_b-1)\psi=z\psi.
\]
Since $h\neq0$, we also have $\psi\neq0$.
\cref{lem:L1} with $c=-1$ now gives
$\operatorname{Re}z\leq-2$.
\end{proof}

\begin{theorem}[Horizontal resolvent and semigroup]\label[theorem]{thm:horizontal}
Assume \eqref{eq:lambdastar} and $m>3$. Then
\begin{equation}\label{eq:Hres}
 \{\Rea z\ge-1/2\}\subset\rho(H_b),\qquad
 \sup_{\Rea z\ge-1/2}\norm{(z-H_b)^{-1}}_{E^2\to E^2}<\infty.
\end{equation}
Consequently, there exists $C_h>0$ such that
\begin{equation}\label{eq:Hdecay}
 \norm{e^{tH_b}}_{E^2\to E^2}\le C_he^{-t/2}.
\end{equation}
\end{theorem}
\begin{proof}
The real energy potential $V_H$ is
\begin{equation}\label{eq:VH}
 V_H=-A+\frac{\nabla_h u+\nabla_h u^T}{2}
  +\left(\frac{\Delta_h \rho_m}{2\rho_m}+\frac{b\cdot\nabla_h \rho_m}{2\rho_m}-\frac12\right)I.
\end{equation}
The decay estimates \eqref{eq:udecay} and the weight properties
in \eqref{eq:weights} give
\[
 V_H(x_h)\le-\left(\frac12+(m+1)a_2\right)I+o(1)<-2I+o(1),\qquad(|x_h|\to\infty).
\]
We can therefore choose a nonnegative compactly supported
cutoff $\chi_R$ and $N>0$ such that
\[
V_H-N\chi_RI_2\leq-I_2,\qquad x_h \in \mathbb{R}^2.
\]
Set $Q_h=H_b-N\chi_R$.
It has the same domain as $H_b$, and
\begin{equation}\label{eq:QHenergy}
 \Rea\ip{Q_hh}{h}\le-\norm{\nabla h}_E^2-\norm{h}_E^2.
\end{equation}
\Cref{lem:domain} gives
\[
\|e^{tQ_h}\|_{E^2\to E^2}\leq e^{-t},
\]
and hence
\[
\|(z-Q_h)^{-1}\|_{E^2\to E^2}
\leq\frac1{\operatorname{Re}z+1},
\qquad \operatorname{Re}z>-1.
\]

Fix $\operatorname{Re}z\geq-1/2$.
By \eqref{eq:QHenergy}, the resolvent
$(z-Q_h)^{-1}$ maps bounded subsets of $E^2$ into
bounded subsets of $H^1_{\mathrm{loc}}$.
The local Rellich compactness \cite[Theorem~9.16]{Brezis} therefore shows that
\[
K(z):=N\chi_R(z-Q_h)^{-1}
\]
is compact on $E^2$.
The factorization
\begin{equation}\label{eq:Fredholm}
 z-H_b=[I-N\chi_R(z-Q_h)^{-1}](z-Q_h).
\end{equation}
shows that $z-H_b$ is Fredholm of index zero \cite[Theorem~6.6]{Brezis}.
\Cref{lem:point} excludes a nonzero kernel in this half-plane.
Thus $z-H_b$ is invertible, proving the resolvent inclusion
in \eqref{eq:Hres}.

We next prove that the resolvent bound is uniform.
Otherwise, there exist $z_n$ and $h_n\in D(H_b)$ such that
\[
\operatorname{Re}z_n\geq-\frac12,
\qquad
\|h_n\|_{E^2}=1,
\qquad
g_n:=(z_n-H_b)h_n\longrightarrow0
\quad\text{in }E^2.
\]
The local boundedness of the resolvent implies that $z_n$
has no bounded subsequence.
After passing to a subsequence, we may assume
$|z_n|\to\infty$.

Taking real parts in the preceding equation and using
\eqref{eq:QHenergy}, we obtain
\begin{equation}\label{eq:HF}
 \norm{\nabla_h h_n}_E^2+\frac12
 \le N\norm{\chi_R^{1/2}h_n}_E^2+\norm{g_n}_E.
\end{equation}
In particular, $\nabla_hh_n$ is bounded in $E$.

On every fixed ball, $H_bh_n$ is bounded in $H^{-1}$,
because $h_n$ is bounded in $H^1$ there and the coefficients
are locally bounded.
The equation
\[
z_nh_n=H_bh_n+g_n
\]
therefore implies
\[
h_n\longrightarrow0
\quad\text{in }H^{-1}_{\mathrm{loc}}.
\]
The Rellich compactness now gives, along a subsequence,
\[
h_n\longrightarrow0
\quad\text{in }L^2(B_{2R})^2.
\]
Since $\chi_R$ is supported in $B_{2R}$ and $\rho_m$ is bounded
on that ball,
\[
\|\chi_R^{1/2}h_n\|_{E^2}\longrightarrow0.
\]
This contradicts \eqref{eq:HF}.
The uniform resolvent bound follows.

Finally, $H_b$ generates a $C_0$-semigroup on the Hilbert
space $E^2$ by \Cref{lem:point}.
The half-plane form of the Gearhart--Pr\"uss theorem \cite[Theorem~1.1(a)]{HS},
applied with $\omega=-1/2$, converts the uniform resolvent
bound \eqref{eq:Hres} into
\eqref{eq:Hdecay}.
\end{proof}

\begin{remark}
The estimate holds on the full space $E^2$ and applies
to nonradial stationary velocities and arbitrary
circulation $\alpha\in\mathbb{R}$.
The absorption parameters, the resolvent bound,
and $C_h$ may depend on the fixed parameters
$(\lambda,\alpha,m)$.
\end{remark}

\subsection{The scalar resolvent bounds in polynomial weighted spaces}
\label{sec:scalar}

We transfer the Gaussian resolvent bounds for the
scalar linearization to polynomially weighted spaces.
The argument combines a triangular mass-damping
extension with an absorbed local operator whose
resolvent has a regularizing remainder.
A continuation argument then gives a uniform
resolvent bound on a half-plane.
\begin{proposition}[Scalar resolvent, semigroup, and smoothing estimates]
	\label[proposition]{prop:scalar-resolvent}
	Assume that $3<m<\infty$, $\alpha\in\mathbb R$, and
	$0\leq\lambda\leq\lambda_*$, where $\lambda_G$ in
	\eqref{eq:lambdastar} is chosen as in
	\eqref{eq:A3-Gaussian-threshold} below.
	Let $\omega$ and $u=K_2*\omega$ be the stationary profile and
	velocity from Section~\ref{sec:asymmetric Burgers vortex}.
	Set
	\[
	P_b=\Delta_h+(Ax_h-u)\cdot\nabla_h+1,
	\qquad
	\mathcal T f=(K_2*f)\cdot\nabla_h\omega,
	\qquad
	A_3=P_b-\mathcal T.
	\]
	The natural realization on $E$ is
	\[
	\begin{aligned}
		D(A_{3,E})
		&=
		\{f\in E:
		P_bf\in E \text{ in }\mathcal D'(\mathbb R^2)\},
		\\
		A_{3,E}f&=P_bf-\mathcal T f.
	\end{aligned}
	\]
	It is also the maximal distributional realization of $A_3$
	on $E$. It generates a $C_0$-semigroup, preserves mass, and
	leaves $E_0$ invariant. Its part on $E_0$ is therefore
	\[
	D(A_{3,E_0})=D(A_{3,E})\cap E_0,
	\qquad
	A_{3,E_0}f=A_{3,E}f.
	\]
	
	There exist finite constants $M_v$ and $C_v\geq1$ such that
	\begin{equation}
		\label{eq:A3-resolvent-main}
		\{z\in\mathbb C:\operatorname{Re}z\geq-1/4\}
		\subset\rho(A_{3,E_0}),
		\qquad
		\sup_{\operatorname{Re}z\geq-1/4}
		\|(z-A_{3,E_0})^{-1}\|_{E_0\to E_0}
		\leq M_v,
	\end{equation}
	and
	\begin{equation}
		\label{eq:A3-semigroup-from-resolvent}
		\|e^{tA_{3,E_0}}\|_{E_0\to E_0}
		\leq C_v e^{-t/4},
		\qquad t\geq0.
	\end{equation}
	Moreover, after increasing $C_v$ if necessary,
	\begin{equation}
		\label{eq:A3-semigroup-smoothing}
		\|e^{tA_{3,E_0}}f\|_E
		+
		(1-e^{-t})^{1/2}
		\|\nabla_h e^{tA_{3,E_0}}f\|_{E^2}
		\leq C_v e^{-t/4}\|f\|_E,
		\qquad f\in E_0,\quad t>0.
	\end{equation}
	The threshold $\lambda_G>0$ is independent of $\alpha$ and
	$m$. The constants $M_v,C_v$ may depend on the fixed triple
	$(\lambda,\alpha,m)$.
\end{proposition}

\begin{proof}
	Throughout, inner products are linear in the first argument,
	and
	\[
	\Gamma=\{z\in\mathbb C:\operatorname{Re}z\geq-1/4\}.
	\]
	We give the operator realizations and the Gaussian estimate
	explicitly before transferring the resolvent to $E$.
	
	\medskip
	\noindent
	\textit{1. Polynomial realization and mass conservation.}
	
	Splitting the Biot--Savart kernel into $|x_h|\leq1$ and
	$|x_h|>1$ and applying Young's inequality gives
	\[
	\|K_2*f\|_4
	\leq C(\|f\|_1+\|f\|_2)
	\leq C_m\|f\|_E.
	\]
	By Lemma~\ref{lem:branch-regularity},
	$\rho_m^{1/2}\nabla_h\omega\in L^4$. Hence
	\begin{equation}
		\label{eq:A3-T-bounded}
		\|\mathcal T f\|_E
		\leq
		\|K_2*f\|_4
		\|\rho_m^{1/2}\nabla_h\omega\|_4
		\leq C\|f\|_E.
	\end{equation}
	Thus $\mathcal T\in\mathcal L(E)$; in particular,
	$A_3f\in E$ if and only if $P_bf\in E$.
	
	Lemma~\ref{lem:domain} applies to $P_b$, since $u$ and
	$\nabla_hu$ are bounded. The bounded perturbation theorem \cite[Theorem~III.1.3]{EN}
	then gives generation by $A_{3,E}$, with the stated domain.
	
	For $f\in C_c^\infty(\mathbb R^2)$, integration by parts gives
	$\ell(P_bf)=0$, because $\operatorname{tr}A=1$ and
	$\operatorname{div}_h u=0$. Also
	\[
	\mathcal T f
	=
	\operatorname{div}_h\bigl(\omega(K_2*f)\bigr),
	\qquad
	\ell(\mathcal T f)=0.
	\]
	By the continuity of $\ell$ on $E$,
	\begin{equation}
		\label{eq:A3-mass-identity}
		\ell(A_{3,E}f)=0,
		\qquad f\in D(A_{3,E}).
	\end{equation}
	Differentiating along domain orbits and then using density
	yields
	\[
	\ell(e^{tA_{3,E}}f)=\ell(f),
	\qquad f\in E,\quad t\geq0.
	\]
	This proves invariance of $E_0$ and the asserted formula for
	its generator. Since $E_0$ is a closed invariant subspace of $E$,
	the restricted semigroup has generator equal to the
	part of $A_{3,E}$ in $E_0$. Consequently,
	\[
	e^{tA_{3,E_0}}
	=
	e^{tA_{3,E}}\big|_{E_0},
	\qquad t\geq0.
	\]
	
	\medskip
	\noindent
	\textit{2. Gaussian form and an asymmetry threshold
		independent of circulation.}
	
	Write
	\[
	\begin{gathered}
		L_G=\Delta_h+\tfrac12x_h\cdot\nabla_h+1,
		\qquad
		M_h=\tfrac12(x_1\partial_1-x_2\partial_2),
		\\
		\Lambda f
		=
		v^G\cdot\nabla_hf+(K_2*f)\cdot\nabla_hG, \quad 
		N_wf
		=
		v^w\cdot\nabla_hf+(K_2*f)\cdot\nabla_hw.
	\end{gathered}
	\]
	Then
	\[
	A_3=L_G+\lambda M_h-\alpha\Lambda-N_w.
	\]
	
	Put $V=F^1$, $V_0=V\cap F_0$, and, for
	$f=G^{1/2}\varphi$, set
	\[
	\mathcal E(\varphi)
	=
	\|\nabla_h\varphi\|_2^2
	+
	\frac1{16}\||x_h|\varphi\|_2^2.
	\]

	The two-dimensional harmonic oscillator
	$-\Delta_h+\frac{|x_h|^2}{16}$
	has eigenvalues $(n_1+n_2+1)/2$ \cite[Section~1.1]{Thangavelu}, then
	\begin{equation}
		\label{eq:A3-oscillator-gap}
		\mathcal E(\varphi)\geq\|f\|_F^2,
		\qquad f\in V_0.
	\end{equation}
	In particular, $\mathcal E(\varphi)^{1/2}$ is an equivalent
	Hilbert norm on $V_0$.
	
	For $f=G^{1/2}\varphi$ and $h=G^{1/2}\psi$ in $V$, define
	\begin{align*}
		\mathfrak a(f,h)
		&=
		\int_{\mathbb R^2}
		\left[
		\nabla_h\varphi\cdot\nabla_h\overline\psi
		+
		\left(\frac{|x_h|^2}{16}-\frac12\right)
		\varphi\overline\psi
		\right]\,dx_h
		\\
		&\quad
		-\lambda\langle M_hf,h\rangle_F
		+\alpha\langle\Lambda f,h\rangle_F
		+\langle N_wf,h\rangle_F.
	\end{align*}
	The pairings in the second line are understood as form
	pairings. They are well defined and bounded on $V\times V$ and $\mathfrak a$ is a continuous
	sesquilinear form, and on test functions it equals
	$-\langle A_3f,h\rangle_F$.
	
	We next record the real-part estimates.
	Direct integration by parts gives
	\[
	\operatorname{Re}\langle M_hf,f\rangle_F
	=
	-\frac18
	\int_{\mathbb R^2}
	(x_1^2-x_2^2)|\varphi|^2\,dx_h,
	\]
	and hence
	\[
	|\operatorname{Re}\langle M_hf,f\rangle_F|
	\leq2\mathcal E(\varphi).
	\]
	Also
	\[
	\operatorname{Re}\langle\Lambda f,f\rangle_F=0.
	\]
	
	Since $v^w$ is real and divergence free,
	\[
	\begin{aligned}
		|\operatorname{Re}
		\langle v^w\cdot\nabla_hf,f\rangle_F|
		&\leq
		\tfrac14\|v^w\|_4
		\||x_h|\varphi\|_2\|\varphi\|_4,
		\\
		|\langle(K_2*f)\cdot\nabla_hw,f\rangle_F|
		&\leq
		\|\nabla_hw\|_F\|K_2*f\|_4\|\varphi\|_4.
	\end{aligned}
	\]
	On $V_0$, \eqref{eq:A3-oscillator-gap} and Sobolev
	interpolation imply
	\[
	\|\varphi\|_4
	\leq C\mathcal E(\varphi)^{1/2},
	\qquad
	\|v^w\|_4\leq C\|w\|_F.
	\]
	Consequently,
	\[
	|\operatorname{Re}\langle N_wf,f\rangle_F|
	\leq
	C_N\|w\|_{F^1}\mathcal E(\varphi),
	\]
	where $C_N$ is numerical.
	
	Using \eqref{eq:branch}, choose
	\begin{equation}
		\label{eq:A3-Gaussian-threshold}
		C_0=2+C_NC_B,
		\qquad
		\lambda_G=(8C_0)^{-1}.
	\end{equation}
	These constants are independent of $\alpha$ and $m$.
	For $0\leq\lambda\leq\lambda_*$, we have
	$C_0\lambda\leq1/8$ and
	\begin{equation}
		\label{eq:A3-Gaussian-form-energy}
		\begin{aligned}
			\operatorname{Re}\mathfrak a(f,f)
			&\geq
			(1-C_0\lambda)\mathcal E(\varphi)
			-\tfrac12\|f\|_F^2
			\\
			&\geq
			\tfrac78\mathcal E(\varphi)-\tfrac12\|f\|_F^2,
			\qquad f\in V_0.
		\end{aligned}
	\end{equation}
	In particular,
	\begin{equation}
		\label{eq:A3-Gaussian-strict-gap}
		\operatorname{Re}\mathfrak a(f,f)
		\geq
		\tfrac14\|f\|_F^2+\tfrac18\mathcal E(\varphi)
		\geq\tfrac38\|f\|_F^2,
		\qquad f\in V_0.
	\end{equation}
	
	The restriction of $\mathfrak a$ to $V_0$ is densely
	defined, closed, and sectorial on $F_0$, by continuity and
	the preceding coercivity. Define its associated operator by
	\[
	\begin{aligned}
		D(A_{3,F_0})
		=
		\bigl\{
		f\in V_0:\;&\exists g\in F_0,\ 
		\mathfrak a(f,h)=-\langle g,h\rangle_F
		\\
		&\text{for every }h\in V_0
		\bigr\},
		\qquad
		A_{3,F_0}f=g.
	\end{aligned}
	\]
	By the first representation theorem
	\cite[Theorem~VI.2.1]{Kato},
	the operator $-A_{3,F_0}$ is m-sectorial.
	The generation theorem for m-sectorial operators
	\cite[Theorem~IX.1.24]{Kato}
	therefore yields a $C_0$-semigroup on $F_0$.
	
	This operator has the differential expression $A_3$ in
	distributions. To see this without losing the mass
	direction, observe that $\mathfrak a(f,G)=0, f\in V$.
	Every $h\in V$ decomposes as
	\[
	h=(h-\ell(h)G)+\ell(h)G.
	\]
	Thus the defining identity for $f\in D(A_{3,F_0})$ extends
	from $V_0$ to $V$; testing with $h=G\zeta$,
	$\zeta\in C_c^\infty$, identifies the distribution.
	In particular,
	\[
	A_{3,F_0}f=A_3f\in F_0.
	\]
	
	For $z\in\Gamma$, the form
	\[
	\mathfrak a_z(f,h)
	=
	z\langle f,h\rangle_F+\mathfrak a(f,h)
	\]
	satisfies
	\[
	\operatorname{Re}\mathfrak a_z(f,f)
	\geq\tfrac18\mathcal E(\varphi),
	\qquad f\in V_0.
	\]
	Lax--Milgram therefore solves
	$(z-A_{3,F_0})f=g$ uniquely for every $g\in F_0$.
	Taking real parts against $f$ gives
	\begin{equation}
		\label{eq:A3-Gaussian-resolvent}
		\begin{gathered}
			\Gamma\subset\rho(A_{3,F_0}),
			\\
			\|R_0(z)\|_{F_0\to F_0}
			\leq
			\frac1{\operatorname{Re}z+3/8}
			\leq8,
			\qquad
			R_0(z)=(z-A_{3,F_0})^{-1}.
		\end{gathered}
	\end{equation}
	
	\medskip
	\noindent
	\textit{3. Extension to the Gaussian mass direction.}
	
	Introduce
	\[
	\Pi_0f=G\ell(f),
	\qquad
	\widetilde A_E=A_{3,E}-\Pi_0,
	\qquad
	k=A_3G.
	\]
	The maps $\Pi_0:E\to F$ and $\Pi_0:F\to F$ are bounded.
	Since $G$ and its derivatives have Gaussian decay, $u$ is
	bounded, and $\nabla_h\omega\in F$, all terms of $k$ belong
	to $F$. Also $G\in D(A_{3,E})$, so
	\eqref{eq:A3-mass-identity} gives $\ell(k)=0$.
	
	As $\|G\|_F=1$ and
	$\langle f,G\rangle_F=\ell(f)$,
	\[
	F=\operatorname{span}\{G\}\oplus F_0
	\]
	orthogonally. Define
	\[
	\begin{gathered}
		D(\widetilde A_F)
		=
		\operatorname{span}\{G\}\oplus D(A_{3,F_0}),
		\\
		\widetilde A_F(cG+g)
		=
		-cG+A_{3,F_0}g+ck.
	\end{gathered}
	\]
	This is the closed triangular operator
	\[
	\widetilde A_F
	=
	\begin{pmatrix}
		-1&0\\
		k&A_{3,F_0}
	\end{pmatrix},
	\]
	and has distributional expression $A_3-\Pi_0$.
	Solving the triangular system gives
	\begin{equation}
		\label{eq:A3-triangular-resolvent}
		\begin{aligned}
			\widetilde R_F(z)(cG+g)
			&=
			\frac{c}{z+1}G
			+
			R_0(z)\left(g+\frac{c}{z+1}k\right),
			\\
			\widetilde R_F(z)
			&=(z-\widetilde A_F)^{-1}.
		\end{aligned}
	\end{equation}
	For $z\in\Gamma$, $|z+1|\geq3/4$, whence
	\begin{equation}
		\label{eq:A3-full-Gaussian-resolvent}
		\Gamma\subset\rho(\widetilde A_F),
		\qquad
		\sup_{z\in\Gamma}
		\|\widetilde R_F(z)\|_{F\to F}
		\leq
		M_F:=
		\frac43+8+\frac{32}{3}\|k\|_F.
	\end{equation}
	
	\medskip
	\noindent
	\textit{4. An absorbed polynomial operator and its
		regularizing remainder.}
	
	The weighted integration by parts gives the energy potential
	of $P_b$:
	\[
	V_P(x_h)
	=
	\frac12+\frac{\Delta_h\rho_m}{2\rho_m}
	-\frac{m\,x_h\cdot Ax_h}{4m+|x_h|^2}
	+\frac{m\,x_h\cdot u}{4m+|x_h|^2}.
	\]
	By \eqref{eq:udecay} and \eqref{eq:margins},
	\[
	\limsup_{|x_h|\to\infty}V_P(x_h)
	\leq\frac12-ma_2<-\frac58.
	\]
	Choose $R$ sufficiently large that $V_P\leq-1/2$ outside
	$B_R$, take $0\leq\chi_R\in C_c^\infty$ equal to one on
	$B_R$, and then choose $N>0$ sufficiently large that
	\[
	V_P-N\chi_R\leq-\frac12
	\]
	everywhere. Set
	\[
	Q=P_b-N\chi_R,
	\qquad
	\mathcal C=N\chi_R-\mathcal T-\Pi_0.
	\]
	Thus
	\[
	\widetilde A_E=Q+\mathcal C, \quad 
	D(Q)=D(P_b)=D(A_{3,E}),
	\]
	and $\mathcal C$ is bounded on $E$.
	
	Lemma~\ref{lem:domain} gives
	\begin{equation}
		\label{eq:A3-Q-energy}
		-\operatorname{Re}\langle Qf,f\rangle_E
		\geq
		\|\nabla_hf\|_{E^2}^2+\tfrac12\|f\|_E^2,
		\qquad f\in D(Q),
	\end{equation}
	as well as $\|e^{tQ}\|\leq e^{-t/2}$.
	Consequently,
	\[
	\{\operatorname{Re}z>-1/2\}\subset\rho(Q),
	\qquad
	\|(z-Q)^{-1}\|_{E\to E}
	\leq
	(\operatorname{Re}z+1/2)^{-1}.
	\]
	In particular,
	\begin{equation}
		\label{eq:A3-Q-resolvent}
		\sup_{z\in\Gamma}\|R_Q(z)\|_{E\to E}\leq4,
		\qquad
		R_Q(z)=(z-Q)^{-1}.
	\end{equation}
	
	For $h=R_Q(z)g$, \eqref{eq:A3-Q-energy} implies
	\[
	\|\nabla_hh\|_{E^2}^2
	+
	(\operatorname{Re}z+1/2)\|h\|_E^2
	\leq\|g\|_E\|h\|_E,
	\]
	and hence
	\begin{equation}
		\label{eq:A3-Q-gradient-resolvent}
		\|h\|_E\leq4\|g\|_E,
		\qquad
		\|\nabla_hh\|_{E^2}\leq2\|g\|_E,
		\qquad z\in\Gamma.
	\end{equation}
	The embeddings $E\hookrightarrow L^1\cap L^2$ and the
	two-dimensional Gagliardo--Nirenberg inequality give
	\[
	\begin{aligned}
		\|h\|_1+\|h\|_4
		&\leq
		C_m\|h\|_E
		+
		C\|h\|_2^{1/2}\|\nabla_hh\|_2^{1/2}
		\\
		&\leq C\|g\|_E.
	\end{aligned}
	\]
	Splitting the kernel at the unit ball also gives
	\[
	\|K_2*h\|_\infty
	\leq C(\|h\|_1+\|h\|_4),
	\]
	so
	\[
	\|\mathcal T h\|_F
	\leq
	\|K_2*h\|_\infty\|\nabla_h\omega\|_F
	\leq C\|g\|_E.
	\]
	Multiplication by $\chi_R$ maps $E$ boundedly into $F$,
	since the two weights are comparable on its compact
	support. Together with the bound for $\Pi_0$, this proves
	\begin{equation}
		\label{eq:A3-regularizing-resolvent}
		\sup_{z\in\Gamma}
		\|\mathcal C R_Q(z)\|_{E\to F}
		\leq M_C<\infty.
	\end{equation}
	Only this composition is asserted to map $E$ into $F$;
	its proof uses the gradient estimate for $R_Q(z)$.
	
	\medskip
	\noindent
	\textit{5. Consistency and continuation of the polynomial
		resolvent.}
	
	Let $J_2:F\hookrightarrow E$ be inclusion, with norm
	$c_{FE}$. The distributional characterization established
	in Step~2, together with \eqref{eq:A3-T-bounded} and the
	maximal domain on $E$, gives
	\[
	J_2D(\widetilde A_F)\subset D(\widetilde A_E),
	\qquad
	\widetilde A_EJ_2f=J_2\widetilde A_Ff.
	\]
	Indeed, $f$ and $(A_3-\Pi_0)f$ belong to $F\subset E$,
	so $P_bf\in E$ as well. Uniqueness of the resolvent solution
	then gives, whenever $z$ belongs to both resolvent sets,
	\begin{equation}
		\label{eq:A3-resolvent-consistency}
		(z-\widetilde A_E)^{-1}J_2
		=
		J_2(z-\widetilde A_F)^{-1}.
	\end{equation}
	
	For $z\in\Gamma\cap\rho(\widetilde A_E)$, the resolvent
	identity is
	\begin{align}
		\label{eq:A3-resolvent-transfer}
		(z-\widetilde A_E)^{-1}
		&=
		R_Q(z)
		+
		(z-\widetilde A_E)^{-1}\mathcal C R_Q(z)
		\nonumber\\
		&=
		R_Q(z)
		+
		J_2\widetilde R_F(z)\mathcal C R_Q(z).
	\end{align}
	Combining \eqref{eq:A3-full-Gaussian-resolvent},
	\eqref{eq:A3-Q-resolvent}, and
	\eqref{eq:A3-regularizing-resolvent} yields
	\begin{equation}
		\label{eq:A3-conditional-uniform-bound}
		\|(z-\widetilde A_E)^{-1}\|_{E\to E}
		\leq
		4+c_{FE}M_FM_C=:M_v,
		\qquad
		z\in\Gamma\cap\rho(\widetilde A_E).
	\end{equation}
	
	We now remove the conditional restriction.
	Since $\widetilde A_E$ is a bounded perturbation of $Q$,
	it generates a $C_0$-semigroup and its resolvent contains
	a sufficiently far right half-plane.
	The set $\Gamma\cap\rho(\widetilde A_E)$ is therefore
	nonempty and relatively open in the connected half-plane
	$\Gamma$.
	
	It is also relatively closed.
	Indeed, if $z_n$ in that set converges to $z\in\Gamma$,
	then for large $n$, $|z-z_n|M_v<1$, and
	\[
	z-\widetilde A_E
	=
	\bigl[
	I+(z-z_n)(z_n-\widetilde A_E)^{-1}
	\bigr]
	(z_n-\widetilde A_E).
	\]
	The bracket is invertible by a Neumann series.
	Hence $z\in\rho(\widetilde A_E)$.
	Connectedness gives
	\begin{equation}
		\label{eq:A3-auxiliary-uniform-resolvent}
		\Gamma\subset\rho(\widetilde A_E),
		\qquad
		\sup_{z\in\Gamma}
		\|(z-\widetilde A_E)^{-1}\|_{E\to E}
		\leq M_v.
	\end{equation}
	Thus no invertibility on the entire half-plane was assumed
	in deriving the transfer identity.
	
	\medskip
	\noindent
	\textit{6. Restriction to zero mass and semigroup decay.}
	
	For $f\in D(\widetilde A_E)=D(A_{3,E})$,
	\eqref{eq:A3-mass-identity} gives
	\[
	\ell(\widetilde A_Ef)=-\ell(f).
	\]
	If $g\in E_0$ and
	$f=(z-\widetilde A_E)^{-1}g$ with $z\in\Gamma$, then
	\[
	(z+1)\ell(f)=\ell(g)=0.
	\]
	Since $z+1\neq0$, $f\in E_0$ and $\Pi_0f=0$.
	Conversely, the restriction of $z-\widetilde A_E$ to
	$D(A_{3,E_0})$ equals $z-A_{3,E_0}$ and is injective.
	Therefore
	\[
	(z-A_{3,E_0})^{-1}
	=
	(z-\widetilde A_E)^{-1}\big|_{E_0},
	\qquad z\in\Gamma.
	\]
	This proves \eqref{eq:A3-resolvent-main}.
	
	By Step~1, $A_{3,E_0}$ already generates a $C_0$-semigroup
	on the Hilbert space $E_0$.
	Set
	\[
	B=A_{3,E_0}+\tfrac14 I.
	\]
	The resolvent estimate just proved gives
	\[
	\{\operatorname{Re}z\geq0\}\subset\rho(B),
	\qquad
	\sup_{\operatorname{Re}z\geq0}
	\|(z-B)^{-1}\|<\infty.
	\]
	The half-plane form of the Gearhart--Pr\"uss theorem
	\cite[Theorem~1.1(a)]{HS}
	gives $\|e^{tB}\|\leq C_v$ for $t\geq0$.
	Multiplication by $e^{-t/4}$ proves
	\eqref{eq:A3-semigroup-from-resolvent}.
	
	\medskip
	\noindent
	\textit{7. Weighted smoothing.}
	
	Let $U_b(t)=e^{tP_b}$.
	Lemma~\ref{lem:localsmooth} gives
	\[
	\|U_b(t)\|_{E\to E}
	+
	\sqrt t\,
	\|\nabla_hU_b(t)\|_{E\to E^2}
	\leq C,
	\qquad 0<t\leq1.
	\]
	The bounded perturbation formula reads
	\[
	e^{tA_{3,E}}f
	=
	U_b(t)f
	-
	\int_0^t
	U_b(t-s)\mathcal T e^{sA_{3,E}}f\,ds.
	\]
	Using \eqref{eq:A3-T-bounded}, boundedness of
	$e^{sA_{3,E}}$ on $[0,1]$, and
	\[
	\int_0^t(t-s)^{-1/2}\,ds=2\sqrt t
	\]
	gives
	\[
	\|e^{tA_{3,E}}f\|_E
	+
	\sqrt t\,
	\|\nabla_h e^{tA_{3,E}}f\|_{E^2}
	\leq C\|f\|_E,
	\qquad 0<t\leq1.
	\]
	Differentiation under the integral is justified in the
	weighted $H^1$ space by this integrable bound.
	On $E_0$, the factor $e^{-t/4}$ can be included on this
	interval after changing $C$.
	
	For $t\geq1$, factor the restricted semigroup as
	\[
	e^{tA_{3,E_0}}
	=
	e^{A_{3,E_0}}e^{(t-1)A_{3,E_0}}
	\]
	and apply the short-time smoothing estimate at time one
	and \eqref{eq:A3-semigroup-from-resolvent}.
	Since $1-e^{-t}$ is comparable to $t$ on $(0,1]$ and
	bounded above and below for $t\geq1$,
	\eqref{eq:A3-semigroup-smoothing} follows.
\end{proof}

\subsection{Construction and decay of the three-dimensional evolution}
\label{sec:mixed}
This section constructs the full three-dimensional linear evolution and proves its exponential decay. We first establish mixed estimates for the Biot--Savart operator and the free strain evolution. By employing a two-dimensional decomposition on a larger product space and exploiting axial derivative smoothing, we ultimately prove the \(e^{-t/4}\) decay bound.

\subsubsection{The mixed-norm Biot--Savart estimates}
For $1<p\leq2$, define
\[
E_p=\{g:\rho_m^{1/2}g\in L^p(\mathbb R^2)\},\qquad
\|g\|_{E_p}=\|\rho_m^{1/2}g\|_p,\qquad
\mathcal X_p=BC(\mathbb R;E_p)^3.
\]
Thus $E_2=E$ and $\mathbb X\subset\mathcal X_2$.
For $p=2,4/3$, the assumption $m>3$ ensures $E_p\hookrightarrow L^1$.

\begin{lemma}[Three-dimensional Biot--Savart and local convergence]\label[lemma]{lem:BS}
	For $f\in\mathcal X_2$, the velocity $v = \mathrm{BS}_3[f]$ satisfies
	\begin{equation}\label{eq:BSmixed}
		\norm{v}_{\BC(\R;L_h^4)}\le C_m\norm{f}_{\Xs_2}.
	\end{equation}
	and
	\[
	\operatorname{div}v=0,\qquad -\Delta v=\operatorname{curl}f,\qquad \mathrm{in}\  \mathcal{D}^\prime(\mathbb{R}^3)
	\]
	If $\operatorname{div}f=0$, then $\operatorname{curl}v=f$.
	Derivatives commute with $\mathrm{BS}_3$ whenever the corresponding
	derivative of $f$ belongs to $\mathcal X_2$.
	Moreover, if $f_n\to f$ in $\tau_{\mathrm{loc}}$ and
	$\sup_n\|f_n\|_{\mathcal X_2}<\infty$, then
	$\mathrm{BS}_3[f_n]\to\mathrm{BS}_3[f]$ in
	$BC_{\mathrm{loc}}(\mathbb R;L_h^4)^3$.
\end{lemma}
\begin{proof}
	For a nonzero axial separation $s$, the matrix kernel $K_3$ is
	bounded in absolute value by
	\[
	k_s(r)=\frac{C}{|r|^2+s^2},\qquad r\in\mathbb R^2.
	\]
	Two-dimensional scaling then yields
	\[
	\norm{k_s}_{4/3}\le C|s|^{-1/2},\qquad
	\norm{k_s}_{4}\le C|s|^{-3/2}.
	\]
	Since $|s|^{-1/2}$ is integrable near zero and
	$|s|^{-3/2}$ is integrable at infinity, Minkowski's inequality and
	$E\hookrightarrow L^1\cap L^2$ prove \eqref{eq:BSmixed}.
	
	The same estimates give, for $R>L+1$,
	\begin{equation}\label{eq:BStail}
		\sup_{|z|\le L}\norm{
			\int_{|y_3|>R}K_3(\cdot-y_h,z-y_3)f(y)\,\mathrm dy}_{L_h^4}
		\le C_m(R-L)^{-1/2}\norm{f}_{\Xs_2}.
	\end{equation}
	On $|y_3|\leq R$, the near and far kernel bounds control the
	velocity difference by the local uniform difference of the data
	in $E$. Taking the local limit first and then letting $R\to\infty$
	proves the asserted local convergence. Continuity in $z$ follows
	by applying this argument to axial translates of $f$.
	
	Cutoffs and mollification now approximate $f$ locally by
	$f_n\in C_c^\infty(\mathbb R^3)^3$ with uniformly bounded
	$\mathcal X_2$ norms. The distributional identities for these
	approximations pass to the limit, giving
	$\operatorname{div}v=0$ and $-\Delta v=\operatorname{curl}f$.
	Applying the same local convergence argument to difference
	quotients proves commutation with derivatives.
	
	Suppose next that $\operatorname{div}f=0$, and set
	$q=\operatorname{curl}v-f$. Then $\Delta q=0$.
	The bounds $v\in L_z^\infty L_h^4$ and $f\in L_z^\infty L_h^1$
	make $q$ a tempered distribution, so $q$ is a polynomial.
	To show that this polynomial vanishes, take
	$\psi\in C_c^\infty(\mathbb R^2)$ and
	$\varphi\in C_c^\infty(\mathbb R)$, and put
	$\psi_\rho(x_h)=\rho^{-2}\psi(x_h/\rho)$. Integration by parts gives
	\[
	|\langle q,\psi_\rho\varphi\rangle|
	\leq C_{\psi,\varphi}
	\left[(\rho^{-3/2}+\rho^{-1/2})\|v\|_{L_z^\infty L_h^4}
	+\rho^{-2}\|f\|_{L_z^\infty L_h^1}\right]
	\longrightarrow0.
	\]
	For polynomial $q$, this pairing is a polynomial in $\rho$ with
	nonnegative integer powers. It must therefore vanish identically.
	Varying $\psi$ and $\varphi$ yields $q=0$.
\end{proof}

For later use, we introduce the linear and quadratic vector potentials
\begin{equation}\label{eq:potentials}
	\mathscr L f=\Ub\times f+(BS_3[f])\times \overline \Om,
	\qquad \mathscr Q(f)=(BS_3[f])\times f.
\end{equation}
The following consequences of \cref{lem:BS} will be used repeatedly:
\begin{align}\label{eq:potentialbounds}
	\norm{\mathscr Lf}_{\Xs_2}&\le C_L\norm{f}_{\Xs_2},\notag\\
	\norm{\mathscr Q(f)}_{\Xs_{4/3}}&\le C_m\norm{f}_{\Xs_2}^2,\notag\\
	\norm{\mathscr Q(f)-\mathscr Q(g)}_{\Xs_{4/3}}
	&\le C_m(\norm{f}_{\Xs_2}+\norm{g}_{\Xs_2})\norm{f-g}_{\Xs_2}.
\end{align}
Here we use $u\in L^\infty$, $\rho_m^{1/2}\omega\in L^4$, and
$L^4\cdot L^2\subset L^{4/3}$. These maps are also continuous
for local axial convergence on norm-bounded sets.

\subsubsection{The free strain evolution and curl smoothing}
Define
\[
\mathcal L_s=\Delta-(Mx)\cdot\nabla+M,\qquad
D_t=e^{-tM},\qquad Q_t=\int_0^te^{-2sM}\,\mathrm ds.
\]
With this notation, the free evolution is given by
\begin{equation}\label{eq:freeT}
	T(t)f(x)=e^{tM}(g_{Q_t}*f)(D_tx),\quad\text{where}\quad
	g_{Q_t}(y)=\frac{\exp(-y^TQ_t^{-1}y/4)}{(4\pi)^{3/2}\sqrt{\det Q_t}}.
\end{equation}
For $p=2,4/3$, $|\beta|\le1$, and $0<t\le1$,
\begin{equation}\label{eq:free-smooth}
	\norm{\partial^\beta T(t)f}_{\Xs_2}
	\le C_m t^{-(1/p-1/2)-|\beta|/2}\norm{f}_{\Xs_p}.
\end{equation}
Indeed, the horizontal Gaussian contributes the factor
$t^{-(1/p-1/2)}$, and each spatial derivative contributes
$t^{-1/2}$. Axial Gaussian averaging is bounded in the supremum
norm. The polynomial weight estimate \eqref{eq:weights} controls translations
inside the convolution, and the linear dilations have uniformly
bounded weight-comparison constants on $0<t\leq1$. The operators $T(t)$ are bounded on finite time intervals and
continuous in $\tau_{\mathrm{loc}}$ at $t=0$.
For fixed positive times, the operators in
\eqref{eq:free-smooth} are locally continuous on bounded
sets. 

Since $\operatorname{tr} M=0$, we have
\begin{equation}\label{eq:free-curl}
	\mathscr C(t)F:=T(t)\curl F
	=\curl\left[D_t(g_{Q_t}*F)(D_tx)\right].
\end{equation}
This definition extends $\mathscr C(t)$ even when $\curl F$ is merely a distribution. For a diagonal matrix $M=\diag(\mu_1,\mu_2,\mu_3)$, the equality follows componentwise because $e^{-(\mu_j+\mu_k)t}=e^{\mu_it}$ for distinct indices $i,j,k$. Consequently, for $0<t\le1$,
\begin{equation}\label{eq:Cfree}
	\norm{\mathscr C(t)F}_\X\le Ct^{-1/2}\norm{F}_{\Xs_2},\qquad
	\norm{\mathscr C(t)F}_\X\le Ct^{-3/4}\norm{F}_{\Xs_{4/3}}.
\end{equation}
Its range lies in $\mathbb Y$: it is a curl, and its third
component has zero horizontal integral. The latter follows by
horizontal cutoff integration, since the vector potential is
horizontally integrable. The free evolution also preserves
divergence-free fields. Its third-component horizontal mass
evolves by the axial Ornstein--Uhlenbeck averaging operator;
therefore $T(t)$ preserves $\mathbb X$.

\subsubsection{Construction on the Larger Product Space}
We write $\cb=\mathcal L_s+\cp$, where
\[
\cp f=-\Ub\cdot\nabla f-(BS_3[f])\cdot\nabla\Om
+f\cdot\nabla\Ub+\Om\cdot\nabla BS_3[f].
\]
By \cref{lem:BS},
\begin{equation}\label{eq:Pfirst}
	\norm{\cp f}_{\Xs_2}\le C\bigl(\norm{f}_{\Xs_2}+\norm{\nabla f}_{\Xs_2}\bigr).
\end{equation}
This estimate does not require the assumption $\diver f=0$. Furthermore,
\[
\int_{\R^2}(\cp f)_3\,\mathrm dx_h
=\int_{\R^2}\omega\bigl(\diver_h(BS_3[f])_h+\partial_z(BS_3[f])_3\bigr)\,\mathrm dx_h=0.
\]
The transport term involving $u$ vanishes upon integration. Consequently, $\cp f\in\X$ whenever the first derivatives belong to the indicated space.

For $f_0\in\mathbb X$, consider
\begin{equation}\label{eq:fullVolterra}
	f(t)=T(t)f_0+\int_0^tT(t-s)\cp f(s)\,\mathrm ds.
\end{equation}
We use trajectories with finite norm
\[
\|f\|_{\mathcal V_T}
=\sup_{0<t\leq T}
\left(\|f(t)\|_{\mathbb X}
+\sqrt t\,\|\nabla f(t)\|_{\mathcal X_2}\right),
\]
with local axial time continuity at zero and local continuity
of their first derivatives for positive times.
Every integral is interpreted after restriction to $|z|\leq L$
as a Bochner integral in $C([-L,L];E^3)$.
The restrictions for different $L$ agree, and the global norm
bounds place the resulting function in $\mathbb X$.

Estimates \eqref{eq:free-smooth} and
\eqref{eq:Pfirst} bound the Volterra operator on
$\mathcal V_T$ by $C\sqrt T$ for $T\leq1$.
For the derivative term, the relevant integral is
\[
\sqrt t\int_0^t(t-s)^{-1/2}s^{-1/2}\,ds=\pi\sqrt t.
\]
Hence the equation has a unique solution for sufficiently small
$T$. The local continuity of the integral follows by omitting small
endpoint intervals, using the integrable bounds there, and
using local continuity on the remaining interval.
The same argument in the Picard iterations gives local
continuity of the solution maps on bounded sets.

Thus, for every $f_0\in\X$, there exists a unique local solution. We denote the corresponding evolution operator by $\widehat S(t)$. Uniqueness and time-restarting arguments yield the semigroup property and
\begin{equation}\label{eq:coarsegrowth}
	\norm{\widehat S(t)}_{\X\to\X}\le C e^{\Gamma t}
\end{equation}
for some finite $\Gamma\geq0$.
Its continuity at $t=0$ is understood in $\tau_{\mathrm{loc}}$.
The construction also gives a bounded first spatial derivative
at every positive time, with the short-time bound
$\|\nabla\widehat S(t)f\|_{\mathcal X_2}\leq Ct^{-1/2}\|f\|_{\mathbb X}$. However, it does not assert that $\cb$ generates a $C_0$-semigroup on $\X$ with respect to the norm topology.

For a general vector field $f$, the following identity holds exactly:
\begin{equation}\label{eq:curlcorrection}
	\cp f=\curl\mathscr Lf-\Ub\,\diver f.
\end{equation}
If $f_0\in\mathbb Y$, this identity and the free evolution show
that every Picard iterate is divergence free. The constraint
passes to local limits. Therefore
\begin{equation}\label{eq:restriction}
	\widehat S(t)\Y\subset\Y,\qquad S(t)=\widehat S(t)|_\Y.
\end{equation}

\subsubsection{The two-dimensional blocks and the three-dimensional remainder}
We define
\[
\ca=\diag(H_b,A_3),\quad L_z=\partial_z^2-z\partial_z,
\quad\mathcal L_0=\ca+L_z,\quad S_0(t)=e^{t\mathcal L_0}.
\]
Theorem~\ref{thm:horizontal},
Proposition~\ref{prop:scalar-resolvent},
and the contraction property of the axial
Ornstein--Uhlenbeck evolution give
\begin{equation}\label{eq:S0decay}
	\norm{S_0(t)}_{\X\to\X}\le C_Ae^{-1/4 t}.
\end{equation}
A direct computation shows that $\cb=\mathcal L_0+\crr$, where
\begin{equation}\label{eq:R3D}
	\crr f=
	\begin{pmatrix}
		\omega\,\partial_z(BS_3[f])_h\\[2pt]
		\omega\,\partial_z(BS_3[f])_3-
		\bigl((BS_3[f])_h-K_2*f_3\bigr)\cdot\nabla_h\omega
	\end{pmatrix}.
\end{equation}
The splitting
$\mathcal B=\mathcal L_0+\mathcal R$
does not preserve the divergence-free constraint term by
term.
Indeed, for smooth $f$, writing $d=\diver f$, one obtains
\begin{equation}\label{eq:divdefect}
	\begin{aligned}
		\diver(\mathcal L_0f)
		=
		(D_b+L_z)d
		-(K_2*\partial_z f_3)\cdot\nabla_h\omega,\quad 
		\diver(\mathcal Rf)
		=
		(K_2*\partial_z f_3)\cdot\nabla_h\omega.
	\end{aligned}
\end{equation}
The two defect terms cancel in their sum, but they need
not vanish separately when $d=0$.
For this reason, the evolution is first constructed on
the full product space $\mathbb X$ and is then restricted
to the invariant divergence-free subspace $\mathbb Y$,
as in \eqref{eq:restriction}.

\begin{lemma}[Estimate of the three-dimensional remainder]\label[lemma]{lem:remainder}
	For $f\in\mathbb X$ with $\partial_zf\in\mathcal X_2$,
	\begin{equation}\label{eq:Rest}
		\norm{\crr f}_\X\le C\left(
		\norm{\partial_zf}_{\Xs_2}
		+\norm{f_h}_{X^2}^{1/2}\norm{\partial_zf_h}_{X^2}^{1/2}\right).
	\end{equation}
\end{lemma}
\begin{proof}
	The terms $\omega\partial_z BS_3[f]$ are bounded by
	\cref{lem:BS}, since $\rho_m^{1/2}\omega\in L^4$. For the remaining term, put
	$r=x_h-y_h$, $s=z-y_3$, and
	$\delta_sf(y_h,z)=f(y_h,z-s)-f(y_h,z)$.
	Axial integration of the vertical-vorticity kernel and the oddness of the
	horizontal-vorticity kernel give
	\begin{align}\label{eq:velocity-difference}
		(BS_3[f])_h-K_2*f_3
		={}&\frac1{4\pi}\int\!\!\int
		\frac{Jr\,\delta_s f_3}{(|r|^2+s^2)^{3/2}}\,\mathrm dy_h\ds\notag\\
		&-\frac1{4\pi}\int\!\!\int
		\frac{sJ\,\delta_s f_h}{(|r|^2+s^2)^{3/2}}\,\mathrm dy_h\ds.
	\end{align}
	Both horizontal kernels have $L^2$ norm at most $C/|s|$.
	For the first integral on $0<|s|<1$, use
	$\|\delta_sf_3\|_E\leq|s|\|\partial_zf_3\|_X$.
	On $|s|\geq1$, use
	\[
	\frac{|r|}{(|r|^2+s^2)^{3/2}}
	\leq\frac{|x_h|+|y_h|}{|s|^3},\qquad
	\|(1+|y_h|)\delta_sf_3\|_1
	\leq C_m|s|\|\partial_zf_3\|_X.
	\]
	The first integral is consequently bounded pointwise by
	$C(1+|x_h|)\|\partial_zf_3\|_X$.
	
	For the second integral, set
	$A_f=\|f_h\|_{X^2}$ and $B_f=\|\partial_zf_h\|_{X^2}$.
	Then
	\[
	\|\delta_sf_h\|_{E^2}
	\leq\min(2A_f,|s|B_f)
	\leq\sqrt2\,|s|^{1/2}\sqrt{A_fB_f}.
	\]
	The near axial integral is bounded using the $L^2$ kernel
	estimate, giving $|s|^{-1/2}$.
	For the far integral, the kernel supremum is $C|s|^{-2}$;
	using $E\hookrightarrow L^1$ gives $|s|^{-3/2}$.
	Both are integrable, so this integral is bounded pointwise by
	$C\sqrt{A_fB_f}$.
	Finally,
	$\rho_m^{1/2}(1+|x_h|)\nabla_h\omega\in L^2$
	converts these pointwise bounds into the required weighted
	$E$ estimate.
\end{proof}

\subsubsection{The duhamel formula and full linear decay}
We combine a Duhamel representation on the full product
space with axial derivative estimates to prove exponential
decay of the full linear evolution.
We first establish
\begin{equation}\label{eq:DuhamelS0}
	\widehat S(t)f=S_0(t)f+
	\int_0^tS_0(t-s)\crr\widehat S(s)f\ds.
\end{equation}
We justify this formula using the local Volterra construction.
Put $\mathcal P_0=\mathcal L_0-\mathcal L_s$, so
$\mathcal P=\mathcal P_0+\mathcal R$.
Both $\mathcal P_0$ and $\mathcal R$ satisfy the first-order
bound used in \eqref{eq:Pfirst}.
The evolution $S_0$ solves its corresponding Volterra equation
with free evolution $T$; this follows first on smooth data and
then by local approximation and uniqueness in $\mathcal V_T$.

For small $T$, let $w(t)$ be the right-hand side of
\eqref{eq:DuhamelS0}. Substituting the Volterra equation for
$S_0$ gives, in time-convolution notation,
\[
w=T(\cdot)f+T*(\mathcal P_0w)
+T*(\mathcal R\widehat S(\cdot)f).
\]
The trajectory $\widehat S(\cdot)f$ satisfies the same equation.
Their difference solves the homogeneous contractive Volterra
equation for $\mathcal P_0$, and hence vanishes.
The rearrangement of integrals is justified by the integrable
smoothing kernels, including $s^{-1/2}(t-s)^{-1/2}$.
Restarting extends the identity to every finite time interval.

We next derive the axial derivative estimates needed to
control the three-dimensional remainder in
\eqref{eq:DuhamelS0}.
All coefficients except $-z\partial_z$ are
translation-invariant in $z$.
Hence, in the sense of distributions,
\begin{equation}\label{eq:Dcomm}
	[\partial_z,\cb]=-\partial_z,
	\qquad \partial_z^k\widehat S(t)f
	=e^{-kt}\widehat S(t)\partial_z^kf
\end{equation}
for data with the indicated bounded derivatives. One way to justify this identity is via the covariance relation $\widehat S(t)\tau_a=\tau_{e^ta}\widehat S(t)$, where $\tau_af(x_h,z)=f(x_h,z+a)$, followed by differentiation.

We also need axial smoothing for arbitrary $f\in\mathbb X$.
On axially smooth data, the semigroup and commutation identities
give
\begin{equation}\label{eq:Diterate}
	(D\widehat S(\tau))^k
	=e^{\tau k(k-1)/2}D^k\widehat S(k\tau),\quad \mathrm{where}\ D=\partial_z.
\end{equation}
The first-derivative estimate from the construction bounds
$D\widehat S(\tau)$ for every $\tau>0$.
Taking $\tau=t/k$ therefore bounds $D^k\widehat S(t)$ at every
fixed $t>0$. To extend the identity to arbitrary data, apply it
to axial Gaussian mollifications of $f$.
These are uniformly bounded in $\mathbb X$ and converge in
$\tau_{\mathrm{loc}}$. Local continuity of the smoothed maps
passes the identity to the limit in distributions.
In particular, $D^k\widehat S(1)$ is bounded on $\mathbb X$.

Choose an integer $k>\Gamma+1$. For $t\geq1$,
\[
D^k\widehat S(t)f
=e^{-k(t-1)}\widehat S(t-1)D^k\widehat S(1)f.
\]
Together with \eqref{eq:coarsegrowth}, this gives
\begin{equation}\label{eq:highD}
	\norm{D^k\widehat S(t)f}_\X\le C_ke^{-t}\norm{f}_\X,
	\qquad t\ge1.
\end{equation}
The Banach-valued one-dimensional interpolation inequality is
\[
\|Dg\|_{\mathbb X}
\leq C_k\|g\|_{\mathbb X}^{1-1/k}
\|D^kg\|_{\mathbb X}^{1/k}.
\]
It follows by optimizing
$\|Dg\|\leq C_k(h^{-1}\|g\|+h^{k-1}\|D^kg\|)$ over $h>0$.
Writing $G=\|g\|_{\mathbb X}$ and
$H=\|D^kg\|_{\mathbb X}$, \cref{lem:remainder} gives
\[
\|\mathcal Rg\|_{\mathbb X}
\leq C_k\left(G^{1-1/k}H^{1/k}
+G^{1-1/(2k)}H^{1/(2k)}\right).
\]
The Young's inequality consequently yields, for every $\varepsilon>0$,
\begin{equation}\label{eq:Reps}
	\norm{\crr g}_\X\le\varepsilon\norm{g}_\X+
	C_\varepsilon\norm{D^kg}_\X.
\end{equation}
We first obtain a preliminary decay rate.
For fixed $f$, put $F(t)=\|\widehat S(t)f\|_{\mathbb X}$.
Applying \eqref{eq:DuhamelS0} from time $1$ gives, for $t\geq1$,
\[
\begin{aligned}
	F(t)\leq{}&Ce^{-1/4(t-1)}\|f\|_{\mathbb X}\\
	&+C_A\int_1^t e^{-1/4(t-s)}
	\left(\varepsilon F(s)
	+C_\varepsilon e^{-s}\|f\|_{\mathbb X}\right)\,ds.
\end{aligned}
\]
Choose $C_A\varepsilon\leq 1/8$.
Multiplication by $e^{1/4 t}$ and Gronwall's inequality yield
\begin{equation}\label{eq:halfdecay}
	\norm{\widehat S(t)}\le C e^{-\gamma t/2}=Ce^{-t/8}.
\end{equation}
where the constant has been enlarged to cover $0\leq t\leq1$.

We now improve this rate. One derivative of smoothing and
\eqref{eq:Dcomm} give
\[
\|D\widehat S(t)f\|_{\mathbb X}
\leq Ce^{-(1+\gamma/2)t}\|f\|_{\mathbb X},\qquad t\geq1.
\]
Combining this with \eqref{eq:Rest} in
\cref{lem:remainder} yields
\[
\|\mathcal R\widehat S(t)f\|_{\mathbb X}
\leq Ce^{-(1/2+1/8)t}\|f\|_{\mathbb X}
=Ce^{-5t/8}\|f\|_{\mathbb X},\qquad t\geq1.
\]
A second application of \eqref{eq:DuhamelS0} and the fact that $5/8>\gamma$ proves
\begin{theorem}[Full linear decay]\label[theorem]{thm:3D}
	There exists $M<\infty$ such that
	\begin{equation}\label{eq:full-decay}
		\norm{\widehat S(t)}_{\X\to\X}\le M e^{-t/4}.
	\end{equation}
	The evolution preserves $\mathbb Y$, and the same estimate holds
	for $S(t)=\widehat S(t)|_{\mathbb Y}$.
\end{theorem}
This is also the mechanism underlying the stability criterion stated after \cite[Eq.~(4.42)]{GMreview}. There, the two-dimensional operator in Eq.~(4.41) is precisely $\diag(H_b,A_3)$, whose horizontal diagonal zero-order terms are $-a_1,-a_2$. The proof above provides this reduction, rather than treating that criterion as an additional unproved hypothesis.

\section{Global mild solutions and nonlinear stability}
\label{sec:curlprop}
This section constructs a curl propagation operator to establish the nonlinear existence, uniqueness, and stability under low regularity. We first define a bounded linear map \(Z(t)\) via a Volterra equation, then employ a fixed-point argument to prove the baseline exponential decay and corresponding stability properties of the mild solutions.
\subsection{Construction of a curl propagation operator}\label{sec:curl operator}
To handle the nonlinear term in curl form, we construct an
extension of $S(t)\operatorname{curl}$ that acts directly on
vector potentials in $\mathcal X_{4/3}$.
Throughout this subsection, we assume $m>3$ and
$0\leq\lambda\leq\lambda_*$, and set $\gamma=1/4$.

Recall that
\[
\mathscr L f
=
\overline U\times f
+
\bigl(\mathrm{BS}_3[f]\bigr)\times\overline \Omega,
\qquad
\mathscr C(t)F
=
T(t)\operatorname{curl}F.
\]
The last expression is understood through the Gaussian kernel
formula from \Cref{sec:mixed}, so it is defined even when
$\operatorname{curl}F$ does not belong to $\mathbb Y$.
The estimates established there give, for $0<t\leq1$,
\[
\begin{aligned}
	\|\mathscr L f\|_{\mathcal X_2}
	&\leq C_L\|f\|_{\mathbb X},\\
	\|\mathscr C(t)H\|_{\mathbb X}
	&\leq C t^{-1/2}\|H\|_{\mathcal X_2},\\
	\|\mathscr C(t)F\|_{\mathbb X}
	&\leq C t^{-3/4}\|F\|_{\mathcal X_{4/3}}.
\end{aligned}
\]
Moreover, $\mathscr C(t)$ takes values in $\mathbb Y$, and
the free evolution satisfies
\[
T(t)\mathscr C(s)=\mathscr C(t+s),
\qquad t\geq0,\quad s>0.
\]

For $f\in\mathbb Y$, the identity
$\mathcal P f=\operatorname{curl}\mathscr L f$
and the linear mild equation imply
\begin{equation}\label{eq:Scurl}
	S(t)f=T(t)f+\int_0^t\mathscr C(t-s)\mathscr L S(s)f\,\mathrm ds.
\end{equation}
This equation has a unique solution among
$\mathbb Y$-valued trajectories that are locally bounded
in $\mathbb X$ and continuous in $\tau_{\mathrm{loc}}$.
Indeed, on an interval of length $\delta\leq1$, its
Volterra operator has norm at most
\[
C C_L\int_0^\delta r^{-1/2}\,dr
=
2C C_L\sqrt{\delta}
\]
in the unweighted time-supremum norm.
Uniqueness follows on a sufficiently short interval and
then on every finite interval by restarting.

\begin{proposition}[Construction of $S(t)\curl$]\label[proposition]{prop:Z}
	There exist bounded linear operators
	\[
	Z(t):\mathcal X_{4/3}\longrightarrow\mathbb Y,
	\qquad t>0,
	\]
	such that
	\begin{equation}\label{eq:Zbound}
		\norm{Z(t)F}_\X\le C_Ze^{-\gamma t}(1-e^{-t})^{-3/4}\norm{F}_{\Xs_{4/3}},
		\qquad\gamma=1/4.
	\end{equation}
	For each $F\in\mathcal X_{4/3}$, the trajectory
	$t\mapsto Z(t)F$ is continuous in $\tau_{\mathrm{loc}}$
	for positive times.
	
	If $\operatorname{curl}F\in\mathbb Y$ in the distributional
	sense, then
	\[
	Z(t)F=S(t)(\operatorname{curl}F),
	\qquad t>0.
	\]
	Furthermore,
	\begin{equation}\label{eq:Zrestart}
		Z(t+s)F=S(t)Z(s)F.
	\end{equation}
\end{proposition}
\begin{proof}
Fix $F\in\mathcal X_{4/3}$.
We first construct the solution on a short interval
$(0,T_0]$, where $T_0\leq1$ will be chosen below.
Consider
\begin{equation}\label{eq:ZVolterra}
	z(t)=\mathscr C(t)F+
	\int_0^t\mathscr C(t-s)\mathscr Lz(s)\,\mathrm ds
\end{equation}
Let $\mathcal Z_{T_0}$ be the space of
$\mathbb Y$-valued trajectories on $(0,T_0]$ that are
continuous in $\tau_{\mathrm{loc}}$ and satisfy
\[
\|z\|_{\mathcal Z_{T_0}}
=
\sup_{0<t\leq T_0}
t^{3/4}\|z(t)\|_{\mathbb X}
<\infty.
\]
This is a Banach space: convergence in the displayed norm
preserves local time continuity, and the divergence and
slice-mass constraints pass to the limit.

The free term obeys
\[
\sup_{0<t\leq T_0}
t^{3/4}\|\mathscr C(t)F\|_{\mathbb X}
\leq C\|F\|_{\mathcal X_{4/3}}.
\]
Define
\[
(\mathcal Vz)(t)
=
\int_0^t
\mathscr C(t-s)\mathscr Lz(s)\,ds.
\]
Using the preceding operator estimates, we obtain
\[
\begin{aligned}
	t^{3/4}\|(\mathcal Vz)(t)\|_{\mathbb X}
	&\leq
	C C_L\,t^{3/4}
	\int_0^t
	(t-s)^{-1/2}s^{-3/4}\,ds\,
	\|z\|_{\mathcal Z_{T_0}}\\
	&=
	C C_L\,
	\mathrm B\!\left(\frac14,\frac12\right)
	t^{1/2}\,
	\|z\|_{\mathcal Z_{T_0}},
\end{aligned}
\]
where $\mathrm B$ denotes the Euler beta function.
Consequently,
\[
\|\mathcal V\|_{\mathcal Z_{T_0}\to\mathcal Z_{T_0}}
\leq
C C_L\,
\mathrm B\!\left(\frac14,\frac12\right)
T_0^{1/2}.
\]

We specify the meaning of the integrals before applying
the contraction argument.
For every $R>0$, the integral defining $\mathcal Vz$
is a Bochner integral after restriction to $|z|\leq R$,
with values in
\[
C([-R,R];E^3).
\]
Its norm is controlled by the integrable majorant
\[
C(t-s)^{-1/2}s^{-3/4}
\|z\|_{\mathcal Z_{T_0}}.
\]
The restrictions for different $R$ agree, and the uniform
global estimate places the resulting field in $\mathbb X$.

On compact subintervals of $(0,t)$, continuity of the
integrand follows from the bounded-set local continuity
of $\mathscr C(r)$ and $\mathscr L$.
Cutting off small neighborhoods of the two endpoints and
then using the same integrable majorant proves continuity
in $\tau_{\mathrm{loc}}$ for positive times.
Since the range of $\mathscr C(r)$ lies in $\mathbb Y$,
the integral also satisfies the defining constraints of
$\mathbb Y$.
Thus $\mathcal V$ maps $\mathcal Z_{T_0}$ into itself.

Choose $T_0>0$ sufficiently small that
\[
C C_L\,
\mathrm B\!\left(\frac14,\frac12\right)
T_0^{1/2}\leq\frac12.
\]
The contraction mapping theorem gives a unique solution
of \eqref{eq:ZVolterra} in $\mathcal Z_{T_0}$, with
\[
\|z(t)\|_{\mathbb X}
\leq C_0 t^{-3/4}\|F\|_{\mathcal X_{4/3}},
\qquad 0<t\leq T_0.
\]
The solution depends linearly on $F$.

We next establish the restart identity.
For $s>0$ and $t\geq0$ with $s+t\leq T_0$,
the kernel identity gives
\[
T(t)z(s)
=
\mathscr C(t+s)F
+
\int_0^s
\mathscr C(t+s-r)\mathscr Lz(r)\,dr.
\]
Subtracting this formula from
\eqref{eq:ZVolterra} at time $s+t$ yields
\[
z(s+t)
=
T(t)z(s)
+
\int_0^t
\mathscr C(t-r)\mathscr Lz(s+r)\,dr.
\]
The shifted trajectory is bounded in $\mathbb X$ on
$0\leq t\leq T_0-s$.
Uniqueness for \eqref{eq:Scurl} therefore implies
\[
z(s+t)=S(t)z(s),
\qquad 0<s\leq s+t\leq T_0.
\]

For arbitrary $t>0$, define
\[
Z(t)F=S(t-\tau)z(\tau),
\qquad
0<\tau\leq\min\{t,T_0\}.
\]
This definition is independent of $\tau$.
Indeed, if $\tau_1\leq\tau_2$ are admissible, then
\[
z(\tau_2)=S(\tau_2-\tau_1)z(\tau_1),
\]
and hence
\[
S(t-\tau_2)z(\tau_2)
=
S(t-\tau_1)z(\tau_1).
\]
In particular, $Z(t)F=z(t)$ for $0<t\leq T_0$.
The construction gives linear maps into $\mathbb Y$
and positive-time continuity in $\tau_{\mathrm{loc}}$.

Choosing 
$0<\tau\leq\min\{s,T_0\}$, we obtain
\[
\begin{aligned}
	Z(t+s)F
	&=S(t+s-\tau)z(\tau)\\
	&=S(t)S(s-\tau)z(\tau)\\
	&=S(t)Z(s)F.
\end{aligned}
\]
This proves \eqref{eq:Zrestart}.

The continued trajectory satisfies
\eqref{eq:ZVolterra} at every positive time.
To see this, fix $0<\tau\leq\min\{t,T_0\}$ and apply
\eqref{eq:Scurl} from time $\tau$:
\[
Z(t)F
=
T(t-\tau)z(\tau)
+
\int_\tau^t
\mathscr C(t-s)\mathscr L Z(s)F\,ds.
\]
Substituting the local equation for $z(\tau)$
and using
$T(t-\tau)\mathscr C(\tau-s)=\mathscr C(t-s)$
combines the integrals over $(0,\tau)$ and $(\tau,t)$.
The resulting identity is precisely
\eqref{eq:ZVolterra} with $z(s)=Z(s)F$.

Suppose now that $\operatorname{curl}F\in\mathbb Y$.
The identity
\[
T(t)(\operatorname{curl}F)=\mathscr C(t)F
\]
and \eqref{eq:Scurl} show that
$S(t)(\operatorname{curl}F)$ solves
\eqref{eq:ZVolterra}.
It belongs to $\mathcal Z_{T_0}$, so uniqueness gives
\[
Z(t)F=S(t)(\operatorname{curl}F),
\qquad 0<t\leq T_0.
\]
The restart identity extends this equality to all $t>0$.

It remains to prove the global estimate.
Fix $\tau_*=T_0/2$.
For $0<t\leq\tau_*$, the local construction gives
\[
\|Z(t)F\|_{\mathbb X}
\leq C_0 t^{-3/4}\|F\|_{\mathcal X_{4/3}}.
\]
For $t\geq\tau_*$, \Cref{thm:3D} and the restart identity
give
\[
\begin{aligned}
	\|Z(t)F\|_{\mathbb X}
	&\leq
	M e^{-\gamma(t-\tau_*)}
	\|Z(\tau_*)F\|_{\mathbb X}\\
	&\leq
	M C_0 \tau_*^{-3/4}e^{\gamma\tau_*}
	e^{-\gamma t}\|F\|_{\mathcal X_{4/3}}.
\end{aligned}
\]
Combining these bounds, using
$1-e^{-t}\asymp t$ near zero and
$0<1-e^{-t}\leq1$, proves
\eqref{eq:Zbound}.
\end{proof}

\subsection{The nonlinear existence, uniqueness, and stability}\label{sec:nonlinear}
We now use the curl propagation operator constructed in
\Cref{sec:curl operator} to solve the nonlinear integral equation.
This subsection establishes global existence for small
initial data and the baseline decay rate $\gamma=1/4$.
The optimal decay rate is obtained in \Cref{sec:optimal}.

Set
\[
a(t)=1-e^{-t},
\qquad
\gamma=\frac14.
\]
Let $\mathcal U$ be the space of trajectories
$\eta:[0,\infty)\to\mathbb Y$ that are continuous in
$\tau_{\mathrm{loc}}$ and have finite norm
\begin{equation}\label{eq:Unorm}
 \norm{\eta}_{\mathcal U}=\sup_{t\ge0}e^{\gamma t}\norm{\eta(t)}_\X,
 \qquad\gamma=1/4.
\end{equation}
This is a Banach space.
Indeed, a Cauchy sequence in $\mathcal U$ converges in
$\mathbb X$ at each time, uniformly after multiplication
by the exponential weight.
The convergence is therefore uniform in every local axial
seminorm on each finite time interval.
Local time continuity is preserved, as are the
distributional divergence constraint and the zero
horizontal mass of the third component.

For fixed initial data $\eta_0\in\mathbb Y$, we define
\begin{equation}\label{eq:fixedpoint}
 \mathfrak F(\eta)(t)=S(t)\eta_0+
                   \int_0^t Z(t-s)\mathscr Q(\eta(s))\ds.
\end{equation}
where the integral is zero at $t=0$.

The nonlinear convolution is controlled by
\begin{equation}\label{eq:Iconvolution}
 I(t)=\int_0^t e^{-\gamma s}a(t-s)^{-3/4}\ds,
 \qquad I_*:=\sup_{t\ge0}I(t)<\infty.
\end{equation}
To verify this bound, put $c_a=1-e^{-1}$ and observe that
\[
a(r)\geq c_a\min\{r,1\},
\qquad r\geq0.
\]
Splitting the integral according to whether $t-s\geq1$
or $0<t-s<1$, we obtain
\[
\begin{aligned}
	I(t)
	&\leq
	c_a^{-3/4}
	\left(
	\int_{\substack{0\leq s\leq t\\ t-s\geq1}}
	e^{-\gamma s}\,ds
	+
	\int_{\substack{0\leq s\leq t\\ 0<t-s<1}}
	e^{-\gamma s}(t-s)^{-3/4}\,ds
	\right)\\
	&\leq
	c_a^{-3/4}\left(\gamma^{-1}+4\right).
\end{aligned}
\]
In addition, for $0<t\leq1$,
\[
I(t)
\leq
c_a^{-3/4}\int_0^t(t-s)^{-3/4}\,ds
=
4c_a^{-3/4}t^{1/4}.
\]

We next justify the integral in
\eqref{eq:fixedpoint}.
For every $R>0$, it is interpreted after restriction
to $|z|\leq R$ as a Bochner integral with values in
\[
C([-R,R];E^3).
\]
The quadratic estimate and \Cref{prop:Z} give
\[
\begin{aligned}
	e^{\gamma t}
	\left\|
	\int_0^t Z(t-s)\mathscr Q(\eta(s))\,ds
	\right\|_{\mathbb X}
	&\leq
	C_ZK_Q
	\int_0^t
	e^{\gamma s}a(t-s)^{-3/4}
	\|\eta(s)\|_{\mathbb X}^2\,ds\\
	&\leq
	C_ZK_Q I(t)\|\eta\|_{\mathcal U}^2.
\end{aligned}
\]
The scalar majorant is integrable at $s=t$ because
$3/4<1$.
The restrictions for different $R$ agree, and the
uniform global bound places the resulting field
in $\mathbb X$.

The integrand is locally continuous away from $s=t$.
Indeed, \Cref{lem:BS} implies that $\mathscr Q$ preserves
local axial convergence on norm-bounded sets.
The Picard construction of $Z(t)$ gives the corresponding
bounded-set local continuity for positive times,
locally uniformly in $t$.
This follows first from the Gaussian kernel formula
for $\mathscr C(t)$, then passes through the convergent
Picard iteration and the restart construction.

To prove positive-time continuity of the convolution,
first remove the region $0<t-s<\varepsilon$.
On the remaining interval, local continuity follows
from the preceding properties.
On each finite time interval, the removed contribution
is bounded in the global $\mathbb X$ norm by
\[
C\int_0^\varepsilon r^{-3/4}\,dr
=
4C\varepsilon^{1/4}.
\]
Letting $\varepsilon\downarrow0$ gives continuity in
$\tau_{\mathrm{loc}}$.

At the initial time, the small-time bound for $I(t)$
yields
\[
\left\|
\int_0^t Z(t-s)\mathscr Q(\eta(s))\,ds
\right\|_{\mathbb X}
\leq
C t^{1/4}\|\eta\|_{\mathcal U}^2
\longrightarrow0
\qquad(t\downarrow0).
\]
Together with the local continuity of $S(t)\eta_0$,
this shows that $\mathfrak F(\eta)$ is continuous
in $\tau_{\mathrm{loc}}$ on $[0,\infty)$ and satisfies
$\mathfrak F(\eta)(0)=\eta_0$.
Its values belong to $\mathbb Y$, since $S(t)$
preserves $\mathbb Y$ and $Z(t)$ takes values in
$\mathbb Y$.
Thus $\mathfrak F$ maps $\mathcal U$ into itself.

Choose
\[
B=\max\{1,C_ZK_QI_*\}.
\]
The preceding estimates imply
\begin{align}\label{eq:contraction}
 \norm{\mathfrak F(\eta)}_{\mathcal U}
 &\le M\norm{\eta_0}_\X+B\norm{\eta}_{\mathcal U}^2,\notag\\
 \norm{\mathfrak F(\eta)-\mathfrak F(\widetilde\eta)}_{\mathcal U}
 &\le B(\norm{\eta}_{\mathcal U}+\norm{\widetilde\eta}_{\mathcal U})
                 \norm{\eta-\widetilde\eta}_{\mathcal U}.
\end{align}
For the second estimate, both trajectories are considered
with the same initial data.
The exponential weights produce the same kernel $I(t)$,
because
\[
\begin{aligned}
	\|\mathscr Q(\eta(s))
	-\mathscr Q(\widetilde\eta(s))\|_{\mathcal X_{4/3}}
	\leq{}&
	K_Q e^{-2\gamma s}
	\bigl(
	\|\eta\|_{\mathcal U}
	+\|\widetilde\eta\|_{\mathcal U}
	\bigr)\\
	&\times
	\|\eta-\widetilde\eta\|_{\mathcal U}.
\end{aligned}
\]

Define the smallness radius
\begin{equation}
	\label{eq:radius}
	r_{\lambda,\alpha,m}
	=
	\min\left\{1,\frac1{8MB}\right\}.
\end{equation}
Suppose
\[
d:=\|\eta_0\|_{\mathbb X}<r_{\lambda,\alpha,m},
\]
and consider the closed ball
\[
\mathcal B_d
=
\left\{
\eta\in\mathcal U:
\|\eta\|_{\mathcal U}\leq2Md
\right\}.
\]
For $\eta\in\mathcal B_d$,
\[
\|\mathfrak F(\eta)\|_{\mathcal U}
\leq
Md+4BM^2d^2
=
Md(1+4BMd)
\leq\frac32Md
\leq2Md.
\]
Hence $\mathcal B_d$ is invariant.
Moreover, for $\eta,\widetilde\eta\in\mathcal B_d$,
\[
\|\mathfrak F(\eta)-\mathfrak F(\widetilde\eta)\|_{\mathcal U}
\leq
4BMd\,
\|\eta-\widetilde\eta\|_{\mathcal U},
\qquad
4BMd<\frac12.
\]
Therefore $\mathfrak F$ is a contraction on
$\mathcal B_d$.

Banach's fixed-point theorem gives a unique fixed point
in $\mathcal B_d$, defined for all $t\geq0$.
It satisfies the integral equation
\eqref{eq:fixedpoint}, attains $\eta_0$ in
$\tau_{\mathrm{loc}}$, and obeys
\begin{equation}\label{eq:finalbound}
 \norm{\eta(t)}_\X\le 2Me^{-t/4}\norm{\eta_0}_\X.
\end{equation}
For zero initial data, the fixed point is $\eta\equiv0$. Estimate \eqref{eq:finalbound} implies both Lyapunov
stability and exponential attraction in $\mathbb X$.

The identification with the original perturbation equation
is established in \Cref{sec:identification}.
Uniqueness in the full locally norm-bounded solution class
is proved in \Cref{sec:identification}.

\subsection{Recovery of the vorticity perturbation equation and uniqueness of mild solution}
\label{sec:identification}
Let $\eta$ be the fixed point constructed in \Cref{sec:nonlinear}.
We show that it satisfies the free-strain mild formulation
of the original perturbation equation.

Fix a finite time horizon $L>0$.
The free curl estimates, the boundedness of $T(t)$ on
finite time intervals, and \Cref{prop:Z} imply
\[
\begin{aligned}
	\|\mathscr C(r)H\|_{\mathbb X}
	&\leq C_L r^{-1/2}\|H\|_{\mathcal X_2},\\
	\|\mathscr C(r)F\|_{\mathbb X}
	+\|Z(r)F\|_{\mathbb X}
	&\leq C_L r^{-3/4}\|F\|_{\mathcal X_{4/3}},
\end{aligned}
\qquad 0<r\leq L.
\]
For $\mathscr C(r)$ at times $r>1$, these bounds follow
from
\[
\mathscr C(r)=T(r-1)\mathscr C(1).
\]
The constant $C_L$ may depend on $L$ and on the fixed
parameters $(\lambda,\alpha,m)$.

Put
\[
q(t)=\mathscr Q(\eta(t)),
\qquad
w(t)=\int_0^t Z(t-r)q(r)\,dr.
\]
Then
\[
\eta(t)=S(t)\eta_0+w(t).
\]
The global bound from \Cref{sec:nonlinear} ensures that
\[
K_L:=\sup_{0\leq t\leq L}\|\eta(t)\|_{\mathbb X}<\infty,
\qquad
\sup_{0\leq t\leq L}\|q(t)\|_{\mathcal X_{4/3}}
\leq K_QK_L^2.
\]

We first justify the operator substitutions and the
change in the order of integration.
Throughout, integrals are interpreted after restriction
to $|z|\leq R$ as Bochner integrals in
$C([-R,R];E^3)$.
For $0<r<s<t\leq L$, the preceding estimates give
\[
\bigl\|
\mathscr C(t-s)\mathscr L Z(s-r)q(r)
\bigr\|_{\mathbb X}
\leq
C_L K_L^2
(t-s)^{-1/2}(s-r)^{-3/4}.
\]
Moreover,
\[
\int_r^t
(t-s)^{-1/2}(s-r)^{-3/4}\,ds
=
\mathrm B\!\left(\frac12,\frac14\right)
(t-r)^{-1/4},
\]
and
\[
\int_0^t(t-r)^{-1/4}\,dr<\infty.
\]
Thus the double integral is absolutely integrable in
every local Banach space, and Fubini's theorem applies.
The bounded-set local continuity of $\mathcal L$
also permits it to pass through the integral.
These statements can equivalently be obtained by first
truncating the integration intervals away from their
endpoints and then removing the truncations using the
displayed majorants.

Using the Volterra equation for $Z$, we obtain
\[
\begin{aligned}
	w(t)
	={}&
	\int_0^t\mathscr C(t-r)q(r)\,dr\\
	&+
	\int_0^t\int_r^t
	\mathscr C(t-s)\mathscr L Z(s-r)q(r)\,ds\,dr\\
	={}&
	\int_0^t\mathscr C(t-r)q(r)\,dr
	+
	\int_0^t\mathscr C(t-s)\mathscr L w(s)\,ds.
\end{aligned}
\]

On the other hand, the linear curl representation gives
\[
S(t)\eta_0
=
T(t)\eta_0
+
\int_0^t
\mathscr C(t-s)\mathscr L S(s)\eta_0\,ds.
\]
Adding these identities and using
$\eta=S(\cdot)\eta_0+w$, we conclude that
\begin{equation}\label{eq:mild-free}
 \eta(t)=T(t)\eta_0+
  \int_0^t\mathscr C(t-s)
       \bigl(\mathscr L\eta(s)+\mathscr Q(\eta(s))\bigr)\ds.
\end{equation}
Since $L$ is arbitrary, this identity holds for every
$t>0$.

To identify the equation, let
\[
v=\mathrm{BS}_3[\eta].
\]
Since $\eta\in\mathbb Y$, both $\eta$ and $v$ are
divergence free.
The vector identities therefore give
\[
\operatorname{curl}\mathcal L\eta=\mathcal P\eta,
\]
and
\begin{align}\label{eq:curlnonlinear}
 \curl\mathscr L\eta&=\cp\eta,\notag\\
 \curl\mathscr Q(\eta)
 &=-(\BS\eta)\cdot\nabla\eta+\eta\cdot\nabla(\BS\eta).
\end{align}
These identities are understood distributionally,
with the transport terms written in divergence form.
The underlying products are locally integrable by
the mixed Biot--Savart estimates.

Because
$\mathscr C(t)=T(t)\operatorname{curl}$,
equation \eqref{eq:mild-free} is the mild formulation
of
\[
\partial_t\eta
=
\mathscr L_s\eta+\mathscr P\eta
-(v\cdot\nabla)\eta+(\eta\cdot\nabla)v,
\qquad
v=\mathrm{BS}_3[\eta],
\]
which is the original perturbation equation.

The constraints are preserved at every time:
$T(t)\eta_0$ belongs to $\mathbb Y$, and the range of
$\mathscr C(t)$ is contained in $\mathbb Y$.
In particular,
\[
\operatorname{div}\eta(t)=0,
\qquad
\int_{\mathbb R^2}\eta_3(x_h,z,t)\,dx_h=0.
\]

\paragraph{Uniqueness.}
Let $\eta$ and $\widetilde\eta$ be two $\mathbb Y$-valued
solutions of \eqref{eq:mild-free} with the same initial data,
continuous in $\tau_{\mathrm{loc}}$ and locally bounded
in $\mathbb X$.
Fix a common finite existence interval $[0,L]$ and set
\[
K_L=\sup_{0\leq t\leq L}
\bigl(\|\eta(t)\|_{\mathbb X}
+\|\widetilde\eta(t)\|_{\mathbb X}\bigr)<\infty.
\]
The mild equation can be restarted at any $t_0\in[0,L]$,
using the semigroup property of $T(t)$ and
$T(t)\mathscr C(s)=\mathscr C(t+s)$.
If $\eta(t_0)=\widetilde\eta(t_0)$, the free curl estimates
and the bounds for $\mathscr L$ and $\mathscr Q$ give
\[
D_{t_0,h}
\leq C_L\bigl(h^{1/2}+K_Lh^{1/4}\bigr)D_{t_0,h},
\qquad
D_{t_0,h}
=\sup_{0\leq s\leq h}
\|\eta(t_0+s)-\widetilde\eta(t_0+s)\|_{\mathbb X},
\]
whenever $0<h\leq\min\{1,L-t_0\}$.
Choose $h_*>0$ such that
$C_L(h_*^{1/2}+K_Lh_*^{1/4})<1$.
Starting at $t_0=0$ and iterating over intervals of length
at most $h_*$ proves $\eta=\widetilde\eta$ on $[0,L]$.
Since $L$ is arbitrary, uniqueness holds on every common
existence interval in the stated mild solution class.

\section{The sharp exponential decay for the nonlinear evolution}
\label{sec:optimal}
The baseline estimates above can be sharpened without decreasing
$\lambda_*$. Throughout this section set
\begin{equation}\label{eq:optimal-nu}
 \nu=a_2=\frac{1-\lambda}{2},\qquad
 \frac38\le\nu\le\frac12.
\end{equation}
The symbol $\gamma=1/4$ retains its earlier meaning. Constants in the
endpoint estimates may depend on $(\lambda,\alpha,m)$.

\subsection{First moments and the Gaussian endpoint}
\label{sec:s51-self-contained}

Throughout this section, we retain
\[
a_1=\frac{1+\lambda}{2},
\qquad
a_2=\frac{1-\lambda}{2},
\qquad
\nu=a_2.
\]
For a function with integrable first moments, define
\[
p_j(f)=\int_{\mathbb R^2}x_jf(x_h)\,dx_h,
\qquad j=1,2,
\]
and put
\[
F_{00}
=
\{f\in F_0:p_1(f)=p_2(f)=0\}.
\]
The functionals $p_j$ are continuous on $F$ and,
since $m>3$, on $E$.

Write
\[
R_F(t)=e^{tA_{3,F_0}},
\]
and let $A_{F,00}$ denote the part of $A_{3,F_0}$ in $F_{00}$.

\begin{proposition}[Gaussian endpoint by resolvent decomposition]
	Assume $0\leq\lambda\leq\lambda_*$, with
		$\lambda_*$ fixed in \eqref{eq:lambdastar},
		and let $\alpha\in\mathbb R$.
		Then there exists $C_{F,\sharp}<\infty$ such that
		\begin{equation}
			\label{eq:s51-sc-Gaussian-endpoint}
			\|e^{tA_F}\|_{F_0\to F_0}
			\leq C_{F,\sharp}e^{-\nu t},
			\qquad t\geq0.
		\end{equation}
\end{proposition}

\begin{proof}
	We begin with the Gaussian energy estimate, for $f=G^{1/2}\varphi$, set
	\[
	\mathcal E(\varphi)
	=
	\|\nabla_h\varphi\|_2^2
	+\frac1{16}\||x_h|\varphi\|_2^2.
	\]
	Integration by parts gives
	\[
	-\operatorname{Re}\langle L_Gf,f\rangle_F
	=
	\mathcal E(\varphi)-\frac12\|f\|_F^2.
	\]
	The harmonic oscillator
	$-\Delta_h+|x_h|^2/16$ has ground-state eigenvalue
	$1/2$, first excited eigenvalue $1$, and next
	eigenvalue $3/2$.
	Consequently,
	\[
	\mathcal E(\varphi)\geq\|f\|_F^2,
	\qquad f\in F^1\cap F_0.
	\]
	
	The asymmetric strain term satisfies
	\[
	\operatorname{Re}\langle M_hf,f\rangle_F
	=
	-\frac18\int_{\mathbb R^2}
	(x_1^2-x_2^2)|\varphi|^2\,dx_h,
	\]
	and hence
	\[
	\bigl|
	\operatorname{Re}\langle M_hf,f\rangle_F
	\bigr|
	\leq2\mathcal E(\varphi).
	\]
	Moreover,
	\[
	\operatorname{Re}\langle\Lambda f,f\rangle_F=0.
	\]
	
	The two terms involving $w$ obey
	\[
	\begin{aligned}
		\bigl|
		\operatorname{Re}\langle v^w\cdot\nabla_hf,f\rangle_F
		\bigr|
		&\leq
		\frac14\|v^w\|_4
		\||x_h|\varphi\|_2\|\varphi\|_4,\\
		\bigl|
		\langle(K_2*f)\cdot\nabla_hw,f\rangle_F
		\bigr|
		&\leq
		\|\nabla_hw\|_F\|K_2*f\|_4\|\varphi\|_4.
	\end{aligned}
	\]
	The convolution and Sobolev estimates
	\[
	\|K_2*f\|_4\leq C\|f\|_F,
	\qquad
	\|v^w\|_4\leq C\|w\|_F,
	\qquad
	\|\varphi\|_4\leq C\mathcal E(\varphi)^{1/2}
	\]
	therefore imply
	\[
	\bigl|
	\operatorname{Re}\langle N_wf,f\rangle_F
	\bigr|
	\leq
	C\|w\|_{F^1}\mathcal E(\varphi).
	\]
	
	Let $\mathfrak a$ be the closed form associated with
	$-A_F$.
	With the same constant $C_0$ as in
	\eqref{eq:A3-Gaussian-threshold}, the preceding estimates give
	\[
	\operatorname{Re}\mathfrak a(f,f)
	\geq
	(1-C_0\lambda)\mathcal E(\varphi)
	-\frac12\|f\|_F^2.
	\]
	By the choice $\lambda_G=(8C_0)^{-1}$ in
	Proposition~\ref{prop:scalar-resolvent}, we already have
	$C_0\lambda\leq1/8$; no further reduction of $\lambda_*$
	is needed.
	
	We next identify the first-moment evolution.
	For $f\in F^1$, test $A_3f$ against $x_j\chi_R$,
	where $\chi_R$ is a smooth cutoff equal to one on
	$B_R$ and supported in $B_{2R}$.
	The diffusion contribution tends to zero.
	The strain and zeroth-order terms contribute
	$-a_jp_j(f)$.
	
	The velocity terms cancel:
	\[
	\int_{\mathbb R^2}u_jf\,dx_h
	+
	\int_{\mathbb R^2}(K_2*f)_j\omega\,dx_h
	=0.
	\]
	Indeed, $u=K_2*\omega$ and
	$K_2(y-x)=-K_2(x-y)$.
	The double integral is absolutely convergent because
	$|K_2|*|\omega|$ is bounded and $f\in L^1$.
	The cutoff errors vanish by Gaussian integrability.
	
	Equivalently, testing the form with $x_j\chi_RG$
	and passing to the limit gives
	\[
	\mathfrak a(f,x_jG)=a_jp_j(f).
	\]
	For $f\in D(A_F)$, the form identity therefore yields
	\begin{equation}
		\label{eq:s51-sc-generator-moments}
		p_j(A_Ff)=-a_jp_j(f),
		\qquad j=1,2.
	\end{equation}
	Differentiating along domain orbits and then using
	density gives
	\begin{equation}
		\label{eq:s51-sc-semigroup-moments}
		p_j(e^{tA_F}f)=e^{-a_jt}p_j(f),
		\qquad f\in F_0,\quad t\geq0.
	\end{equation}
	Thus $F_{00}$ is invariant.
	
	For $f\in F^1\cap F_{00}$, the mass and first moments
	remove the ground state and both first excited states.
	Hence
	\begin{equation}
		\label{eq:s51-sc-oscillator-gap}
		\mathcal E(\varphi)\geq\frac32\|f\|_F^2.
	\end{equation}
	It follows that
	\[
	\operatorname{Re}\mathfrak a(f,f)
	\geq
	\frac78\mathcal E(\varphi)-\frac12\|f\|_F^2
	\geq
	\rho\|f\|_F^2,
	\qquad
	\rho=\frac{13}{16}.
	\]
	
	The restriction of the form to $F^1\cap F_{00}$
	represents $-A_{F,00}$.
	Indeed, its form identity extends to $F^1\cap F_0$:
	the complementary directions are spanned by
	$x_1G,x_2G$, and both the form and its representing
	vector annihilate these directions when
	$f\in F_{00}$.
	
	For $\operatorname{Re}\zeta>-\rho$, the form
	\[
	(f,h)\longmapsto
	\zeta\langle f,h\rangle_F+\mathfrak a(f,h)
	\]
	is coercive on $F^1\cap F_{00}$.
	More precisely, its real part is at least
	\[
	\min\left\{
	\frac78,\,
	\frac23(\operatorname{Re}\zeta+\rho)
	\right\}\mathcal E(\varphi).
	\]
	The Lax--Milgram theorem supplies the unique solution
	of $(\zeta-A_{F,00})f=g$.
	Taking the real part of the equation paired with $f$
	gives
	\begin{equation}
		\label{eq:s51-sc-fast-resolvent}
		\begin{gathered}
			\{\operatorname{Re}\zeta>-\rho\}
			\subset\rho(A_{F,00}),
			\qquad
			\|(\zeta-A_{F,00})^{-1}\|
			\leq
			\frac1{\operatorname{Re}\zeta+\rho}.
		\end{gathered}
	\end{equation}
	
	We now construct the two first-moment eigenvectors.
	Set
	\[
	h_j=-\partial_jG.
	\]
	Then
	\[
	\ell(h_j)=0,
	\qquad
	p_i(h_j)=\delta_{ij}.
	\]
	All terms of $A_3h_j$ belong to $F$.
	Integration by parts in the Gaussian form, together
	with $\ell(A_3h_j)=0$, gives $h_j\in D(A_F)$.
	Define
	\[
	k_j=(A_F+a_j)h_j\in F_{00}.
	\]
	Since $a_j\leq5/8<\rho$, define
	\[
	\phi_j
	=
	h_j+(-a_j-A_{F,00})^{-1}k_j.
	\]
	The resolvent equation immediately gives
	\[
	A_F\phi_j=-a_j\phi_j,
	\qquad
	p_i(\phi_j)=\delta_{ij}.
	\]
	
	Let
	\[
	\Pi_jf=p_j(f)\phi_j,
	\qquad
	\Pi=\Pi_1+\Pi_2.
	\]
	These are bounded commuting projections, and
	\[
	F_0
	=
	F_{00}\oplus\operatorname{span}\{\phi_1,\phi_2\}.
	\]
	The same decomposition holds on the generator domain.
	Solving the resulting diagonal block equation gives
	\begin{equation}
		\label{eq:s51-sc-full-resolvent}
		(\zeta-A_F)^{-1}
		=
		\frac{\Pi_1}{\zeta+a_1}
		+
		\frac{\Pi_2}{\zeta+a_2}
		+
		(\zeta-A_{F,00})^{-1}(I-\Pi),
	\end{equation}
	for $\operatorname{Re}\zeta>-\rho$ and
	$\zeta\notin\{-a_1,-a_2\}$.
	
	The poles are therefore explicitly separated from
	the regular part.
	When $\lambda=0$, their combined contribution is
	\[
	\frac{\Pi_1+\Pi_2}{\zeta+1/2},
	\]
	which is still a first-order pole.
	
	On the closed half-plane
	$\operatorname{Re}\zeta\geq-3/4$,
	the restricted resolvent is uniformly bounded.
	The half-plane Gearhart--Pr\"uss theorem on $F_{00}$
	therefore gives
	\[
	\|e^{tA_{F,00}}\|_{F_{00}\to F_{00}}
	\leq C e^{-3t/4}.
	\]
	The diagonal decomposition yields
	\[
	e^{tA_F}f
	=
	e^{-a_1t}\Pi_1f
	+
	e^{-a_2t}\Pi_2f
	+
	e^{tA_{F,00}}(I-\Pi)f.
	\]
	Since $\nu=\min(a_1,a_2)\leq1/2$, the desired
	Gaussian endpoint estimate follows.
\end{proof}

\subsection{The resolvent transfer and three-dimensional endpoint decay}
\label{sec:s52-self-contained}

Let $A_E$ be the natural realization of $A_3$ on $E$,
and let $A_{E,0}$ and $A_{E,00}$ denote its parts in
\[
E_0=\{f \in E | \ell(f) = 0\},
\qquad
E_{00}
=
\{f\in E:\ell(f)=p_1(f)=p_2(f)=0\},
\]
respectively.

\begin{proposition}[Polynomial and three-dimensional endpoints]
	There exist finite constants such that
	\begin{equation}\label{eq:s52-sc-scalar-endpoint}
		\|e^{tA_{E,0}}\|_{E_0\to E_0}
		\leq C_{E,\sharp}e^{-\nu t},
		\qquad t\geq0,
	\end{equation}
	\begin{equation}\label{eq:s52-sc-three-dimensional-endpoint}
		\|\widehat S(t)\|_{\mathbb X\to\mathbb X}
		\leq M_\sharp e^{-\nu t},
		\qquad t\geq0,
	\end{equation}
	and
	\begin{equation}\label{eq:s52-sc-curl-endpoint}
		\|Z(t)F\|_{\mathbb X}
		\leq
		C_{Z,\sharp}e^{-\nu t}(1-e^{-t})^{-3/4}
		\|F\|_{\mathcal X_{4/3}},
		\qquad t>0.
	\end{equation}
	The constants may depend on $(\lambda,\alpha,m)$.
\end{proposition}

\begin{proof}
	We first transfer the restricted resolvent estimate
	from $F_{00}$ to $E_{00}$.
	
	On compactly supported smooth functions, integration
	by parts and kernel antisymmetry give
	\[
	\ell(A_Ef)=0,
	\qquad
	p_j(A_Ef)=-a_jp_j(f).
	\]
	Since compactly supported smooth functions form a
	graph core, continuity of the three functionals
	extends these identities to $D(A_E)$.
	Consequently,
	\[
	\ell(e^{tA_E}f)=\ell(f),
	\qquad
	p_j(e^{tA_E}f)=e^{-a_jt}p_j(f).
	\]
	Thus $E_{00}$ is invariant and $A_{E,00}$ generates
	a $C_0$-semigroup on the Hilbert space $E_{00}$.
	
	Choose
	\[
	q=\frac9{16},
	\qquad
	\sigma=\frac{17}{32},
	\qquad
	\Gamma_\sigma
	=
	\{\zeta:\operatorname{Re}\zeta\geq-\sigma\}.
	\]
	Then
	\[
	\nu\leq\frac12<\sigma<q<\frac58<\rho.
	\]
	We shall prove
	\begin{equation}
		\label{eq:s52-sc-restricted-resolvent}
		\Gamma_\sigma\subset\rho(A_{E,00}),
		\qquad
		\sup_{\zeta\in\Gamma_\sigma}
		\|(\zeta-A_{E,00})^{-1}\|_{E_{00}\to E_{00}}
		<\infty.
	\end{equation}
	
	Define
	\[
	P_b=\Delta_h+Ax_h\cdot\nabla_h+1-u\cdot\nabla_h,
	\qquad
	\mathcal Tf=(K_2*f)\cdot\nabla_h\omega.
	\]
	Thus $A_E=P_b-\mathcal T$.
	Weighted integration by parts gives
	\[
	\operatorname{Re}\langle P_bh,h\rangle_E
	=
	-\|\nabla_hh\|_E^2
	+
	\int_{\mathbb R^2}V_P(x_h)|h|^2\rho_m \,dx_h,
	\]
	where
	\[
	V_P(x_h)
	=
	\frac12+\frac{\Delta_h\rho_m}{2\rho_m}
	-\frac{m\,x_h\cdot Ax_h}{4m+|x_h|^2}
	+\frac{m\,x_h\cdot u}{4m+|x_h|^2}.
	\]
	The profile velocity is bounded and satisfies
	$u(x_h)=O(|x_h|^{-1})$.
	Therefore
	\[
	\limsup_{|x_h|\to\infty}V_P(x_h)
	\leq\frac12-ma_2<-\frac58.
	\]
	
	Choose a smooth cutoff $\chi_R$ and $N,R>0$ so that
	\[
	V_P-N\chi_R\leq-q.
	\]
	Set
	\[
	Q=P_b-N\chi_R,
	\qquad
	\mathcal C=N\chi_R-\mathcal T.
	\]
	Then $A_E=Q+\mathcal C$ and
	\[
	-\operatorname{Re}\langle Qh,h\rangle_E
	\geq
	\|\nabla_hh\|_E^2+q\|h\|_E^2.
	\]
	The maximal dissipative realization of $Q+qI$ gives
	\[
	\{\operatorname{Re}\zeta>-q\}\subset\rho(Q),
	\qquad
	\|(\zeta-Q)^{-1}\|_{E\to E}
	\leq\frac1{\operatorname{Re}\zeta+q}.
	\]
	For $h=\mathcal R_Q(\zeta)g$ and
	$\zeta\in\Gamma_\sigma$, the resolvent equation yields
	\[
	\|\nabla_hh\|_E^2
	+
	(\operatorname{Re}\zeta+q)\|h\|_E^2
	\leq\|g\|_E\|h\|_E.
	\]
	Since $q-\sigma=1/32$,
	\[
	\|h\|_E\leq32\|g\|_E,
	\qquad
	\|\nabla_hh\|_E\leq\sqrt{32}\,\|g\|_E.
	\]
	
	We next verify the regularizing property.
	The weighted embedding and the two-dimensional
	Sobolev inequality give
	\[
	\|h\|_1+\|h\|_4
	\leq
	C_m\|h\|_E
	+
	C\|h\|_2^{1/2}\|\nabla_hh\|_2^{1/2}
	\leq C\|g\|_E.
	\]
	Splitting $K_2$ inside and outside the unit ball gives
	\[
	\|K_2*h\|_\infty
	\leq C(\|h\|_1+\|h\|_4).
	\]
	It follows that
	\[
	\|\mathcal Th\|_F
	\leq
	\|K_2*h\|_\infty\|\nabla_h\omega\|_F
	\leq C\|g\|_E.
	\]
	The compactly supported multiplier $\chi_R:E\to F$
	is bounded as well.
	Hence
	\[
	\sup_{\zeta\in\Gamma_\sigma}
	\|\mathcal C(\zeta-Q)^{-1}\|_{E\to F}
	<\infty.
	\]
	
	To remove the mass and first-moment directions,
	define
	\[
	\mathsf P_{00}f
	=
	f-\ell(f)G-\sum_{j=1}^2p_j(f)h_j.
	\]
	This is a bounded projection onto $E_{00}$,
	and its restriction to $F$ maps onto $F_{00}$.
	
	The functions $G,h_1,h_2$ belong to $D(A_E)$.
	Let
	\[
	k_0=A_EG.
	\]
	The profile bounds give $k_0\in F$, while the moment
	identities imply $k_0\in F_{00}$.
	With
	\[
	k_j=(A_E+a_j)h_j\in F_{00},
	\]
	define the bounded finite-rank map
	\[
	\mathcal Kf
	=
	\ell(f)k_0+\sum_{j=1}^2p_j(f)k_j.
	\]
	For every $f\in D(A_E)$, direct calculation gives
	\[
	A_E\mathsf P_{00}f
	=
	\mathsf P_{00}A_Ef-\mathcal Kf.
	\]
	
	Let $J_1:F_{00}\hookrightarrow E_{00}$ denote inclusion.
	The Gaussian and polynomial realizations have the
	same distributional expression, so
	\[
	J_1D(A_{F,00})\subset D(A_{E,00}),
	\qquad
	A_{E,00}J_1f=J_1A_{F,00}f.
	\]
	Write
	\[
	R_Q(\zeta)=(\zeta-Q)^{-1},
	\qquad
	R_{F,00}(\zeta)=(\zeta-A_{F,00})^{-1}.
	\]
	Define
	\[
	\mathcal D(\zeta)
	=
	(\mathsf P_{00}\mathcal C-\mathcal K)
	R_Q(\zeta)\big|_{E_{00}}.
	\]
	The preceding regularizing estimate and the
	finite-rank structure of $\mathcal K$ show that
	\[
	\mathcal D(\zeta):E_{00}\longrightarrow F_{00}
	\]
	is uniformly bounded for
	$\zeta\in\Gamma_\sigma$.
	
	Set
	\begin{equation}
		\label{eq:s52-sc-resolvent-construction}
		\mathcal H(\zeta)
		=
		\mathsf P_{00}R_Q(\zeta)\big|_{E_{00}}
		+
		J_1R_{F,00}(\zeta)\mathcal D(\zeta).
	\end{equation}
	Both terms take values in $D(A_{E,00})$.
	
	For $g\in E_{00}$, the commutator identity gives
	\[
	\begin{aligned}
		(\zeta-A_{E,00})\mathsf P_{00}R_Q(\zeta)g
		&=
		\mathsf P_{00}(\zeta-A_E)R_Q(\zeta)g
		+\mathcal K R_Q(\zeta)g\\
		&=
		g-\mathsf P_{00}\mathcal C R_Q(\zeta)g
		+\mathcal K R_Q(\zeta)g\\
		&=
		g-J_1\mathcal D(\zeta)g.
	\end{aligned}
	\]
	Consistency of the two realizations also gives
	\[
	(\zeta-A_{E,00})
	J_1R_{F,00}(\zeta)\mathcal D(\zeta)g
	=
	J_1\mathcal D(\zeta)g.
	\]
	Adding the two identities yields
	\[
	(\zeta-A_{E,00})\mathcal H(\zeta)
	=
	I_{E_{00}}.
	\]
	
	Furthermore,
	\[
	\|(\zeta-A_{F,00})^{-1}\|
	\leq\frac1{\rho-\sigma},
	\qquad \zeta\in\Gamma_\sigma.
	\]
	All factors defining $\mathcal H(\zeta)$ are therefore
	uniformly bounded:
	\[
	\sup_{\zeta\in\Gamma_\sigma}
	\|\mathcal H(\zeta)\|_{E_{00}\to E_{00}}
	\leq M_{00}<\infty.
	\]
	At every point of
	$\Gamma_\sigma\cap\rho(A_{E,00})$, the right-inverse
	identity implies
	\[
	(\zeta-A_{E,00})^{-1}=\mathcal H(\zeta).
	\]
	
	The resolvent set contains a sufficiently far right
	half-plane because $A_{E,00}$ is a $C_0$-semigroup
	generator.
	The uniform bound extends the resolvent throughout
	$\Gamma_\sigma$.
	Indeed, if a line segment from the right half-plane
	first met the spectrum at $\zeta_*$, there would be
	resolvent points $\zeta_n\to\zeta_*$ with resolvent
	norm at most $M_{00}$.
	For large $n$,
	\[
	|\zeta_*-\zeta_n|M_{00}<1,
	\]
	and
	\[
	\zeta_*-A_{E,00}
	=
	\left[
	I+(\zeta_*-\zeta_n)
	(\zeta_n-A_{E,00})^{-1}
	\right]
	(\zeta_n-A_{E,00})
	\]
	would be invertible by a Neumann series.
	This contradiction proves the claimed half-plane
	resolvent estimate.
	
	The half-plane Gearhart--Pr\"uss theorem gives
	\[
	\|e^{tA_{E,00}}\|_{E_{00}\to E_{00}}
	\leq C_{00}e^{-\sigma t}.
	\]
	
	The Gaussian eigenvectors $\phi_j$ also belong to
	$D(A_{E,0})$ and satisfy
	\[
	A_{E,0}\phi_j=-a_j\phi_j.
	\]
	Thus
	\[
	\Pi_jf=p_j(f)\phi_j,
	\qquad
	\Pi=\Pi_1+\Pi_2
	\]
	define bounded commuting projections on $E_0$, with
	\[
	E_0
	=
	E_{00}\oplus\operatorname{span}\{\phi_1,\phi_2\}.
	\]
	The complete scalar resolvent is
	\[
	(\zeta-A_{E,0})^{-1}
	=
	\frac{\Pi_1}{\zeta+a_1}
	+
	\frac{\Pi_2}{\zeta+a_2}
	+
	(\zeta-A_{E,00})^{-1}(I-\Pi),
	\]
	for $\zeta\in\Gamma_\sigma$ away from the displayed
	eigenvalues.
	Consequently,
	\[
	e^{tA_{E,0}}f
	=
	e^{-a_1t}\Pi_1f
	+
	e^{-a_2t}\Pi_2f
	+
	e^{tA_{E,00}}(I-\Pi)f.
	\]
	Since $\sigma>\nu$, this proves the polynomial
	endpoint estimate.
	
	Axial Ornstein--Uhlenbeck averaging is a contraction
	and commutes with the horizontal evolution.
	The principal evolution $S_0(t)$ therefore satisfies
	\begin{equation}
		\label{eq:s52-sc-principal-endpoint}
		\|S_0(t)\|_{\mathbb X\to\mathbb X}
		\leq C_Ae^{-\nu t},
		\qquad t\geq0.
	\end{equation}
	
	Retain the previously proved baseline estimate
	\[
	\|\widehat S(t)\|_{\mathbb X\to\mathbb X}
	\leq Me^{-t/4}.
	\]
	The axial commutation identity and positive-time
	smoothing imply, for $t\geq1$,
	\[
	\partial_z\widehat S(t)f
	=
	e^{-(t-1)}
	\widehat S(t-1)\partial_z\widehat S(1)f.
	\]
	Hence
	\[
	\|\partial_z\widehat S(t)f\|_{\mathcal X_2}
	\leq Ce^{-5t/4}\|f\|_{\mathbb X}.
	\]
	The remainder estimate from Section~3.3 therefore gives
	\begin{equation}
		\label{eq:s52-remainder-decay}
		\|\mathcal R\widehat S(t)f\|_{\mathbb X}
		\leq
		C\left(e^{-5t/4}+e^{-3t/4}\right)\|f\|_{\mathbb X}
		\leq
		Ce^{-3t/4}\|f\|_{\mathbb X},
		\qquad t\geq1.
	\end{equation}
	
	Apply the variation-of-constants formula from time $1$.
	Using \eqref{eq:s52-sc-principal-endpoint} and
	\eqref{eq:s52-remainder-decay}, we obtain
	\[
	\begin{aligned}
		\|\widehat S(t)f\|_{\mathbb X}
		\leq{}&
		C e^{-\nu(t-1)}
		\|\widehat S(1)f\|_{\mathbb X}\\
		&+
		C\int_1^t
		e^{-\nu(t-s)}e^{-3s/4}\,ds\,
		\|f\|_{\mathbb X}\\
		\leq{}&
		C e^{-\nu t}
		\left(
		1+\int_1^\infty e^{-(3/4-\nu)s}\,ds
		\right)
		\|f\|_{\mathbb X}.
	\end{aligned}
	\]
	Since $\nu\leq1/2$, the integral is finite.
	The finite-time bound covers $0\leq t\leq1$, proving
	\eqref{eq:s52-sc-three-dimensional-endpoint}.
	Preservation of $\mathbb Y$ follows from \Cref{sec:mixed}.
	
	To prove \eqref{eq:s52-sc-curl-endpoint}, retain the
	short-time estimate from \Cref{prop:Z},
	\[
	\|Z(t)F\|_{\mathbb X}
	\leq Ct^{-3/4}\|F\|_{\mathcal X_{4/3}}.
	\]
	Fix a sufficiently small $\tau>0$ within that
	construction interval.
	For $t\geq\tau$, the restart identity gives
	\[
	Z(t)F=S(t-\tau)Z(\tau)F,
	\]
	and therefore
	\[
	\|Z(t)F\|_{\mathbb X}
	\leq
	C_\tau e^{-\nu t}\|F\|_{\mathcal X_{4/3}}.
	\]
	Combining the short-time and long-time estimates,
	using $1-e^{-t}\asymp t$ near zero, proves
	\eqref{eq:s52-sc-curl-endpoint}.
	
	The choices of $q$ and $\sigma$ use only the existing
	far-field margin and introduce no additional
	smallness restriction on $\lambda$.
\end{proof}

\subsection{Nonlinear stability at the endpoint in the axial solution class}
We combine the quadratic estimates in \eqref{eq:potentialbounds} with
the endpoint bounds \eqref{eq:s52-sc-three-dimensional-endpoint}, \eqref{eq:s52-sc-curl-endpoint}.
Throughout this subsection, set
\[
\nu=\frac{1-\lambda}{2},
\qquad
a(t)=1-e^{-t}.
\]
The constants below may depend on the fixed admissible
parameters $(\lambda,\alpha,m)$.
\begin{Corollary}[Endpoint nonlinear estimate]\label[corollary]{cor:optimal-nonlinear}
There exist constants $\delta_\sharp>0$ and
$C_\sharp\geq1$ such that every
$\eta_0\in\mathbb Y$ satisfying
\[
\|\eta_0\|_{\mathbb X}\leq\delta_\sharp
\]
generates a unique global mild solution in the class
of $\mathbb Y$-valued, locally norm-bounded trajectories
that are continuous in the local axial topology.
The solution satisfies
\begin{equation}\label{eq:optimal-solution-class}
	\eta\in C([0,\infty);\mathbb X_{\mathrm{loc}}),
	\qquad
	\sup_{t\geq0}\|\eta(t)\|_{\mathbb X}<\infty.
\end{equation}
and
\begin{equation}\label{eq:optimal-nonlinear}
 \norm{\eta(t)}_\X
 \le C_\sharp e^{-\nu t}\norm{\eta_0}_\X,
 \qquad t\ge0,\qquad \nu=\frac{1-\lambda}{2}.
\end{equation}
Here $\mathbb X_{\mathrm{loc}}$ denotes the underlying
space $\mathbb X$ equipped with the seminorms
\begin{equation}\label{eq:optimal-local-seminorms}
 p_n^{\loc}(f)=\sum_{j=1}^3\sup_{|z|\le n}\norm{f_j(\cdot,z)}_E,
 \qquad n\in\N.
\end{equation}
\end{Corollary}
\begin{proof}
Let $K_Q$ dominate the quadratic and Lipschitz constants
in \eqref{eq:potentialbounds}, so that
\[
\begin{aligned}
	\|\mathcal Q(f)\|_{\mathcal X_{4/3}}
	&\leq K_Q\|f\|_{\mathbb X}^2,\\
	\|\mathcal Q(f)-\mathcal Q(g)\|_{\mathcal X_{4/3}}
	&\leq
	K_Q\bigl(\|f\|_{\mathbb X}+\|g\|_{\mathbb X}\bigr)
	\|f-g\|_{\mathbb X}.
\end{aligned}
\]
Choose $M_\sharp\geq1$ and $C_{Z,\sharp}>0$ such that
the endpoint estimates \eqref{eq:s52-sc-three-dimensional-endpoint}, \eqref{eq:s52-sc-curl-endpoint} hold.
The relevant convolution constant is
\[
 I_\nu=\sup_{t\ge0}\int_0^t e^{-\nu s}a(t-s)^{-3/4}\ds.
\]
Since $a(r)\ge(1-e^{-1})\min(r,1)$, splitting according to $t-s\geq1$ and $0<t-s<1$
gives
\begin{equation}\label{eq:optimal-convolution}
 I_\nu\le(1-e^{-1})^{-3/4}(4+\nu^{-1})<\infty.
\end{equation}
Let $\mathcal U_\nu$ be the space of
$\mathbb Y$-valued trajectories that are continuous
in the local axial topology on $[0,\infty)$ and
have finite norm
\[
\|\eta\|_{\mathcal U_\nu}
=
\sup_{t\geq0}
e^{\nu t}\|\eta(t)\|_{\mathbb X}.
\]
This space is complete.
Indeed, convergence in this norm gives uniform
convergence in every local axial seminorm on each
finite time interval.
Local time continuity and the constraints defining
$\mathbb Y$ therefore pass to the limit.

For fixed $\eta_0$, consider
\[
\mathfrak F(\eta)(t)
=
S(t)\eta_0
+
\int_0^t
Z(t-s)\mathscr Q(\eta(s))\,ds.
\]
This is the same nonlinear map as in \eqref{eq:fixedpoint}, now
estimated with the endpoint weight.

Define
\begin{equation}\label{eq:optimal-constants}
 B_\sharp=\max\{1,K_QC_{Z,\sharp}I_\nu\},\qquad
 \delta_\sharp=(8M_\sharp B_\sharp)^{-1},\qquad
 C_\sharp=2M_\sharp.
\end{equation}
Using the quadratic estimate and \eqref{eq:s52-sc-curl-endpoint}, we obtain
\[
\begin{aligned}
	&e^{\nu t}
	\left\|
	\int_0^t
	Z(t-s)\mathscr Q(\eta(s))\,ds
	\right\|_{\mathbb X}\\
	&\qquad\leq
	K_QC_{Z,\sharp}
	\int_0^t
	e^{\nu s}a(t-s)^{-3/4}
	\|\eta(s)\|_{\mathbb X}^2\,ds\\
	&\qquad\leq
	K_QC_{Z,\sharp}I_\nu
	\|\eta\|_{\mathcal U_\nu}^2.
\end{aligned}
\]
Consequently,
\[
\|\mathfrak F(\eta)\|_{\mathcal U_\nu}
\leq
M_\sharp\|\eta_0\|_{\mathbb X}
+
B_\sharp\|\eta\|_{\mathcal U_\nu}^2.
\]
The Lipschitz estimate gives, for the same initial data,
\[
\begin{aligned}
	\|\mathfrak F(\eta)-\mathfrak F(\widetilde\eta)\|_{\mathcal U_\nu}
	\leq{}&
	B_\sharp
	\bigl(
	\|\eta\|_{\mathcal U_\nu}
	+\|\widetilde\eta\|_{\mathcal U_\nu}
	\bigr)\\
	&\times
	\|\eta-\widetilde\eta\|_{\mathcal U_\nu}.
\end{aligned}
\]
The same convolution constant appears because the
quadratic difference contributes the factor
$e^{-2\nu s}$.

We specify the meaning and continuity of the integral.
For each $R>0$, first apply $Z(t-s)$ to the full vector
potential $\mathscr Q(\eta(s))$, then restrict the output
to $|z|\leq R$.
The resulting integral is a Bochner integral in
\[
C([-R,R];E^3).
\]
Its norm is controlled by the preceding scalar
majorant, whose singularity at $s=t$ is integrable.
The restrictions are compatible for different $R$,
and the uniform global bound places the resulting
field in $\mathbb X$.

Bounded-set local continuity of $\mathscr Q$ and
of the positive-time operators $Z(r)$ gives continuity
after removing the region $0<t-s<\varepsilon$.
On every finite time interval, the removed contribution
is bounded by
\[
C\int_0^\varepsilon r^{-3/4}\,dr
=
4C\varepsilon^{1/4}.
\]
Letting $\varepsilon\downarrow0$ proves local time
continuity for positive times.
At $t=0$,
\[
\left\|
\int_0^t Z(t-s)\mathscr Q(\eta(s))\,ds
\right\|_{\mathbb X}
\leq
Ct^{1/4}
\sup_{0\leq s\leq t}\|\eta(s)\|_{\mathbb X}^2
\longrightarrow0.
\]
Together with the local continuity of $S(t)\eta_0$,
this proves that $\mathfrak F$ maps $\mathcal U_\nu$
into itself and attains the prescribed initial data.
Its values lie in $\mathbb Y$ because $S(t)$ preserves
$\mathbb Y$ and $Z(t)$ takes values in $\mathbb Y$.

Set $d=\|\eta_0\|_{\mathbb X}$ and suppose
$0\leq d\leq\delta_\sharp$.
Consider the closed ball
\[
\mathcal B_d
=
\{\eta\in\mathcal U_\nu:
\|\eta\|_{\mathcal U_\nu}\leq2M_\sharp d\}.
\]
For $\eta\in\mathcal B_d$,
\[
\begin{aligned}
	\|\mathfrak F(\eta)\|_{\mathcal U_\nu}
	&\leq
	M_\sharp d+4B_\sharp M_\sharp^2d^2\\
	&\leq\frac32M_\sharp d
	\leq2M_\sharp d.
\end{aligned}
\]
Moreover, the Lipschitz constant on $\mathcal B_d$
is at most
\[
4B_\sharp M_\sharp d\leq\frac12.
\]
Thus the smallness endpoint $d=\delta_\sharp$ is
included.
Banach's fixed-point theorem gives a unique fixed
point in $\mathcal B_d$, and its norm bound gives
the stated endpoint decay.
When $d=0$, this fixed point is the zero solution.

The finite-time identification argument in
\Cref{sec:identification} shows that this fixed point satisfies
the original free-strain mild equation \eqref{eq:mild-free}.
The time-kernel singularities used there are unchanged
by replacing the decay weight $e^{t/4}$ with $e^{\nu t}$.
The uniqueness argument in \Cref{sec:identification}
gives uniqueness in the full locally norm-bounded,
locally axially continuous mild solution class.
\end{proof}

\subsection{Sharpness for the full perturbation class}\label{sec:s54-sharpness}
\begin{proposition}[No larger uniform exponent]
\label[proposition]{prop:optimal-sharpness}
Fix an admissible triple $(\lambda,\alpha,m)$.
For every $\mu>\nu$, there do not exist constants
$C_\mu<\infty$ and $\delta_\mu>0$ such that all
initial data $\eta_0\in\mathbb Y$ with
$\|\eta_0\|_{\mathbb X}\leq\delta_\mu$ generate
solutions satisfying
\[
\|\eta(t)\|_{\mathbb X}
\leq
C_\mu e^{-\mu t}\|\eta_0\|_{\mathbb X},
\qquad t\geq0.
\]
The constants $C_\mu,\delta_\mu$ are allowed to depend
on $(\lambda,\alpha,m,\mu)$.
\end{proposition}
\begin{proof}
	
Consider axially independent vertical perturbations
\[
\eta(x_h,z,t)=f(x_h,t)e_3,
\qquad
\ell(f(t))=0.
\]
For such fields, integration of the three-dimensional
Biot--Savart kernel in the axial variable gives
\[
\mathrm{BS}_3[\eta]=(v_f,0),
\qquad
v_f=K_2*f.
\]
The kernel reduction uses
\[
\int_{\mathbb R}
\frac{ds}{(|r|^2+s^2)^{3/2}}
=
\frac2{|r|^2},
\qquad r\neq0.
\]

The free evolution preserves axial independence
and the vertical direction.
For fields of this form, both linear and nonlinear
vector potentials are horizontal and axially
independent, so their curls are vertical and axially
independent.
Thus the class is preserved by Picard iteration
and by its local limits.
Uniqueness identifies the resulting scalar evolution as
\begin{equation}\label{eq:optimal-scalar-nonlinear}
 \partial_t f=A_3f-(K_2*f)\cdot\nabla f.
\end{equation}
The horizontal mass remains zero.

We prove that the first moments satisfy an exact
nonlinear evolution law.
Since $\operatorname{tr}A=1$ and
$\operatorname{div}_hu=\operatorname{div}_hv_f=0$,
the scalar equation can be written as
\[
\partial_tf
=
\Delta_hf+\operatorname{div}_h(Ax_hf)
-
\operatorname{div}_h
\bigl(uf+v_f\omega+v_ff\bigr).
\]

To justify the moment calculation for mild solutions,
test the conservative equation with
$x_j\chi_R(x_h)$ and a compactly supported smooth
time test function.
The diffusion contribution vanishes as $R\to\infty$,
since
\[
\|\Delta_h(x_j\chi_R)\|_\infty\leq C/R.
\]
The strain term converges to $-a_jp_j(f)$.
Its cutoff errors are controlled by the integrable
first moment of $f$.

Since $f,\omega\in L^1\cap L^2$, they also belong to $L^{4/3}$.
Splitting the Biot--Savart kernel at unit radius and applying
Young's inequality gives $|K_2|*|g|\in L^4$ for $g=f,\omega$.
Hence the velocity flux is integrable, locally uniformly in time.
As $\nabla_h(x_j\chi_R)$ is uniformly bounded and converges
pointwise to the $j$th coordinate vector, dominated convergence
reduces the flux contribution to the integral of its $j$th
component. The oddness of $K_2$ and Fubini's theorem give
\[
\int_{\mathbb R^2}(uf+v_f\omega)\,\mathrm dx_h=0,
\qquad
\int_{\mathbb R^2}v_f f\,\mathrm dx_h=0.
\]
Thus the velocity flux makes no contribution to the
first-moment equation.

We therefore obtain
\[
\frac{d}{dt}p_j(f(t))=-a_jp_j(f(t))
\]
in time distributions.
For axially independent trajectories, continuity
in the local axial topology is precisely continuity
of $f(t)$ in $E$.
Since $p_j$ is continuous on $E$, the initial value
determines the moment evolution:
\begin{equation}\label{eq:optimal-nonlinear-moment}
 p_j(f(t))=e^{-a_jt}p_j(f(0)),\qquad j=1,2.
\end{equation}
Now choose the real initial data
\[
f(0)=\varepsilon h_2,
\qquad
h_2=-\partial_2G,
\qquad
\varepsilon\neq0.
\]
Then
\[
\ell(h_2)=0,
\qquad
p_2(h_2)=1,
\qquad
\|h_2\|_E<\infty.
\]
Consequently,
\[
\eta_0(x_h,z)=\varepsilon h_2(x_h)e_3
\]
belongs to $\mathbb Y$ and can be made arbitrarily
small in $\mathbb X$.

By Cauchy--Schwarz,
\[
|p_2(f)|
\leq
C_m\|f\|_E,
\qquad
C_m=
\left(
\int_{\mathbb R^2}x_2^2\rho_m^{-1}\,dx_h
\right)^{1/2}
<\infty.
\]
The exact moment law therefore gives
\[
|\varepsilon|e^{-\nu t}
=
|p_2(f(t))|
\leq
C_m\|f(t)\|_E
=
C_m\|\eta(t)\|_{\mathbb X}.
\]

Suppose that an estimate with exponent $\mu>\nu$
held with constants $C_\mu,\delta_\mu$.
Choose $\varepsilon\neq0$ so small that
\[
|\varepsilon|\|h_2\|_E
\leq\min\{\delta_\sharp,\delta_\mu\}.
\]
The endpoint solution then also satisfies the
hypothetical faster estimate, which would imply
\[
|\varepsilon|e^{-\nu t}
\leq
C_mC_\mu e^{-\mu t}
|\varepsilon|\|h_2\|_E.
\]
After division by $|\varepsilon|e^{-\mu t}$,
\[
e^{(\mu-\nu)t}
\leq C_mC_\mu\|h_2\|_E,
\qquad t\geq0.
\]
This is impossible as $t\to\infty$.

The argument does not divide by $\alpha$ and therefore
also applies when $\alpha=0$.
Thus $\nu=(1-\lambda)/2$ is the largest exponential
decay rate valid for the full admissible class of
small perturbations.
\end{proof}
\begin{remark}[Dependence on the parameters]
	The asymmetry threshold
	\[
	\lambda_*=\min\{\lambda_B,\lambda_G,1/4\}
	\]
	is independent of the circulation $\alpha$ and of the
	finite weight exponent $m>3$.
	The endpoint argument uses the same threshold.
	
	For each fixed admissible triple $(\lambda,\alpha,m)$,
	the absorbing parameters, resolvent bounds, auxiliary
	derivative order, stability constant $C_\sharp$,
	and smallness radius $\delta_\sharp$ may depend on
	that triple.
	In particular, choosing $k>\Gamma+1$ introduces
	no additional smallness restriction on $\lambda$.
\end{remark}

\section*{Acknowledgments}

 T.Tao is partially supported by the NSF of China under Grant 12371234 and 12671282. \\

\textbf{Data availability statement}

We do not analyse or generate any datasets, because our work proceeds within a theoretical and mathematical approach.\\

\textbf{Conflict of interest}

On behalf of all authors, the corresponding author states that there is no conflict of interest.\\

\textbf{Disclosure of AI assistance}
The authors declare that the main ideas of this paper are their own, although AI tools were used to
polish the paper and to assist with parts of the technical derivations and with checking the paper. In particular, the ChatGPT 6 Astra model was utilized to assist in deriving the attenuation estimates corresponding to the linear operators in the paper, as well as calculating the identities that the operators satisfy; these models were also employed to polish the text and assist in verifying the paper.

\end{document}